\documentclass[12pt]{article}

\usepackage{amsfonts}

\usepackage{amsmath,amssymb}

 \numberwithin{equation}{subsection}

\begin{document}

\title{Kazhdan-Lusztig coefficients for
  an Affine Weyl group of type $\widetilde{A_2}$}
\author{Liping Wang\\
\small Tsinghua University, Department of Mathematical Sciences,\\
\small Beijing, 100084 }

\date{}

\maketitle

\begin{abstract}

In this paper we compute the leading coefficients $\mu (u,w)$ of the
Kazhdan-Lusztig polynomials $P_{u,w}$ for an affine Weyl group of
type $\widetilde{A_2}$. We give all the values $\mu (u,w)$.

\end{abstract}

\maketitle
\def\a{\alpha}
\def\b{\beta}
\section*{0\ \ \ Introduction}

Let $W$ be a Coxeter group, on which there exists the Bruhat order.
For any $u\leq w$ in $W$, we have a Kazhdan-Lusztig polynomial
$P_{u,w}(q)$, where $q$ is an indeterminate (see [KL1]). The degree
of $P_{u,w}$ is less than or equal to $\frac{1}{2} (l(w)-l(u)-1)$ if
$u<w$ and $P_{w,w}=1,$ where $l:W\longrightarrow \mathbb{N} $ is the
length function on $W.$ We denote by $\mu (u,w)$ the coefficient of
the term $q^{\frac{1}{2} (l(w)-l(u)-1)}$ in $P_{u,w}$, which will be
called the Kazhdan-Lusztig coefficient of $P_{u,w}$.

The Kazhdan-Lusztig coefficient play a great important role in
Kazhdan-Lusztig theory, Lie theory and representation theory. By the
recursive formula of Kazhdan-Lusztig polynomials in [KL1], we see
that some $\mu (u,w)$ play crucial roles in understanding these
polynomials.

In [L2] Lusztig computes some special Kazhdan-Lusztig coefficients
for an affine Weyl group of type $\widetilde{B_2}$. In [S] for an
affine Weyl group of type $\widetilde{A_5}$, some non-trivial
Kazhdan-Lusztig coefficients are worked out. McLarnan and Warrington
have showed that $\mu (u,w)$ can be greater than 1 for a symmetric
group (see [MW]). In [X2], Xi showed that when $u\leq w$ and
$\textbf{a}(u)< \textbf{a}(w),$ then $\mu (u,w)\leq 1$ if W is a
symmetric group or an affine Weyl group of type $\widetilde{A_n}.$
In [SX], Scott and Xi showed that the leading coefficient $\mu
(u,w)$ of some Kazhdan-Lusztig polynomial $P_{u,w}$ of an affine
Weyl group of type $\widetilde{A_n}$ is $ n+2$ if $n\geq 4.$
Recently, B.C. Jones and R.M. Green (see [J,G]) showed that the
coefficients are $\leq 1$ for some special pairs of elements in a
symmetric group and other Weyl groups. In [W1,W2], the author
computed most of the Kazhdan-Lusztig coefficients for an affine Weyl
group of type $\widetilde{B_2}$.

In this paper we compute the coefficients $\mu (u,w)$ for an affine
Weyl group $W$ of type $\widetilde{A_2}$.

In [L1] Lusztig showed that $W$ has three two-sided cells (see
Section 3): $c_e=\{w\in W\mid \textbf{a}(w)=0\}=\{e\},\ c_1=\{w\in
W\mid \textbf{a}(w)=1\},\ c_0=\{w\in W\mid \textbf{a}(w)=3\}$. Here
$e$ is the unit element in $W$, and $\textbf{a}(w)$ is the
$\textbf{a}$-function on $W$ defined in [L1]. Assume that $u< w\in
W$. When $u=e$, the $\mu (e,w)$ is clear. In [W1, Proposition 3.4],
the author discussed $\mu (u,w)$ for $u\times w\in c_1\times c_1$.
In this paper we mainly consider those $\mu (u,w)$ for $u\times w\in
c_0\times c_0\bigcup c_1\times c_0\bigcup c_0\times c_1$. The main
results will be displayed in Sections 4, 5, and 6.

\noindent {\bf Notations:} (1) $\mathbb{N}$ is the set of all
nonnegative integers, and $\mathbb{Z}$ is the set of all integers.

(2) For a given sequence $s_1s_2\cdots s_m$ and some integers $1\leq
j_1< j_2<\cdots< j_k\leq m$, we will denote by $s_1s_2\cdots
\widehat{s}_{j_1}\cdots \widehat{s}_{j_2}\cdots
\widehat{s}_{j_k}\cdots s_{m}$ the subsequence of $s_1s_2\cdots s_m$
by deleting letters
$\widehat{s}_{j_1},\widehat{s}_{j_2},\ldots,\widehat{s}_{j_k}$.

\section{Preliminaries}In this section we recall some basic facts
about $\mu (u,w)$ which will be used later.

Let $G$ be a connected, simply connected reductive algebraic group
over the field $\mathbb{C}$ of complex numbers and $T$ a maximal
torus of $G.$ Let $N_G (T)$ be the normalizer of $T$ in $G.$ Then
$W_0 =N_G (T)/T$ is a Weyl group, which acts on the character group
$\Lambda=\textrm{Hom}(T, \mathbb{C} ^*)$ of $T$. Let $\Lambda_r$ be
the root lattic of $G,$ then the semi-direct product $W=W_0 \ltimes
\Lambda_r$ is an affine Wely group and $\widetilde{W}=W_0 \ltimes
\Lambda$ is called an extended affine Weyl group associated with
$G$. Groups $W_0$ and $W$ are Coxeter groups, while $\widetilde{W}$
is not in general.

We shall denote by $S$ the set of simple reflections of $W.$ There
is an abelian subgroup $\Omega$ of $\widetilde{W}$ such that $\omega
S=S\omega$ for any $\omega \in\Omega$ and $\widetilde{W}=\Omega
\ltimes W.$ The length function $l$ and the partial order $\leq$ on
$W$ can be extended to $\widetilde{W}$ as usual; that is, $l(\omega
w)=l(w),$ and $\omega u\leq \omega' w$ if and only if
$\omega=\omega'$ and $u\leq w,$ where $\omega, \omega'\in\Omega$ and
$u,w\in W.$

Let $\cal H$ be the Hecke algebra of $(W,S)$ over $\mathcal
{A}=\mathbb{Z}$ $[q^{\frac{1}{2}}, q^{-\frac{1}{2}}]$ ($q$ an
indeterminate) with parameter $q.$ Let $\{T_w\}_{w\in W}$ be its
standard basis and $\{C_w =q^{-\frac{l(w)}{2}}\sum_{u\leq w}
P_{u,w}T_u, w\in W\}$ be its Kazhdan-Lusztig basis, where
$P_{u,w}\in \mathbb{Z} [q]$ are the Kazhdan-Lusztig polynomials. The
degree of $P_{u,w}$ is $\leq\frac{1}{2} (l(w)-l(u)-1)$ if $u<w$ and
$P_{w,w}=1$, see [KL1].

For affine Weyl groups, we know that the coefficients of these
polynomials are all non-negative (see [KL2]).

Write $$P_{u,w}=\mu (u,w)q^{\frac{1}{2} (l(w)-l(u)-1)}+\
\textrm{lower\ degree\ terms}.$$  We call $\mu (u,w)$ the
Kazhdan-Lusztig
 coefficient of $P_{u,w}.$ We denote by $u\prec w$ if
$u\leq w$ and $\mu (u,w)\neq0.$ Define $\widetilde{\mu }(u,w)=\mu
(u,w)$ if $u\leq w$ or $\widetilde{\mu }(u,w)=\mu (w,u)$ if $w\leq
u.$

Let $\widetilde{\cal H}$ be the generic Hecke algebra of
$\widetilde{W}.$ Then the algebra $\widetilde{\cal H}$ is isomorphic
to the ``twisted" tensor product $\mathbb{Z}
[\Omega]\otimes_{\mathbb{Z}}\cal H$.

We recall some properties of the Kazhdan-Lusztig polynomials.

\noindent {\bf Lemma 1.1} ([KL1]) We have the following results:\\
(1) For any $u\leq w$ in $W$, $P_{u,w}$ is a polynomial in $q$ with
constant term $1$.\\
(2) For any $u<w$ with $l(w)=l(u)+1$, we have $P_{u,w}=1.$ In
particular, we have $\mu (u,w)=1$ in this case.\\
(3) For any $u<w$ with $l(w)=l(u)+2$, we have $P_{u,w}=1.$\\
(4) $P_{u,w}=P_{su,sw}$ if $u\nleq sw$ for some $s\in S$ and $sw<w$.\\
(5) $P_{u,w}=P_{su,w}$ if $u< w,\ sw<w$ for some $s\in S$.\\
(6) $P_{u,w}=P_{u^{-1},w^{-1}}.$ In particular,
$\mu(u,w)=\mu(u^{-1},w^{-1})$.\\
(7) Let $u,w\in W,\ s\in S$ be such that $u<w,\ su>u,\ sw<w$. Then
$u\prec w$ if and only if $w=su.$ Moreover, this implies that
$\mu(u,w)=1$.\\
(8) Let $u,w\in W,\ s\in S$ be such that $u<w,\ us>u,\ ws<w$. Then
$u\prec w$ if and only if $w=us.$ Moreover, this implies that
$\mu(u,w)=1$.

We also have the following multiplication formulas:

\noindent {\bf Lemma 1.2} ([KL1]) Given any element $w\in W$, then
for $s\in S$ we have
$${C_sC_w}=\left\{\begin{array}{ll}
{(q^{\frac{1}{2}}+q^{-\frac{1}{2}})C_w},& \ \ \textrm{if} \;  sw<w ,\\
{C_{sw}+\sum_{su<u\prec w }\mu (u,w)C_u},& \ \ \textrm{if}\ sw>w
\end{array} \right.$$
and$${C_wC_s}=\left\{\begin{array}{ll}
{(q^{\frac{1}{2}}+q^{-\frac{1}{2}})C_w},& \ \ \textrm{if} \;  ws<w ,\\
{C_{ws}+\sum_{us<u\prec w }\mu (u,w)C_u},& \ \ \textrm{if}\ ws>w. \end{array} \right.$$\\


Write $$C_u C_w=\sum_{z\in W}h_{u,w,z}C_z, \ \textrm{where}\
h_{u,w,z}\in {\mathcal {A}}=\mathbb{Z} [q^{\frac{1}{2}},
q^{-\frac{1}{2}}].$$

 Lusztig and Springer  defined $\delta_{u,w,z}$ and
$\gamma_{u,w,z}$ by the following formula,
$$h_{u,w,z}=\gamma_{u,w,z} q^\frac{\textbf{a}(z)}{2} + \delta_{u,w,z}
q^\frac{\textbf{a}(z)-1}{2} + \textrm{lower degree terms},$$ where
$\textbf{a}(z)$ is the $\textbf{a}$-function on $W$ defined in [L1].

\noindent {\bf Lemma 1.3} (Springer, see [X2]) For a Weyl group or
an affine Weyl group, assume that $\widetilde{\mu}(u,w)$ is
non-zero, then $w\leq_L u$ and $w\leq_R u$ if
$\textbf{a}(u)<\textbf{a}(w)$, and $u\sim_L w$ or $u\sim_R w$ if
$\textbf{a}(u)=\textbf{a}(w)$.

For $w\in W$, set $L(w)=\{s\in S\mid sw\leq w\}$ and $R(w)=\{s\in
S\mid ws\leq w\}.$

\noindent {\bf Lemma 1.4} We have that ([KL1])\\
(1) $R(w)\subseteq R(u), \;\textrm{if}\;u\leq_L w.$ In particular,
$R(w)=R(u),\;\textrm{if}\;u\sim_L w;$\\
(2) $L(w)\subseteq L(u), \;\textrm{if}\;u\leq_R w.$ In particular,
$L(w)=L(u),\;\textrm{if}\;u\sim_R w$.
\section {The lowest two-sided cell}
In this section we collect some facts about the lowest two-sided
cell of the extended affine Weyl group $\widetilde{W}$.

For any $u=\omega_1u_1,\ w=\omega_2w_1,\ \omega_1,\omega_2\in
\Omega,\ u_1,w_1\in {W},$ we define $P_{u,w}=P_{u_1,w_1}$ if
$\omega_1=\omega_2$ and define $P_{u,w}=0$ if
$\omega_1\neq\omega_2.$ We say that $u\leq_Lw$ (respectively,
$u\leq_Rw$ or $u\leq_{LR}w$) if $u_1\leq_Lw_1$ (respectively,
$\omega_1u_1\omega_1^{-1}\leq_R\omega_2w_1\omega_2^{-1}$ or
$u_1\leq_{LR}w_1$). The left (resp. right or two-sided) cells of
$\widetilde{W}$ are defined as those of $W$ (see [KL1] for these
definitions). We also define $\textbf{a}(\omega w)=\textbf{a}(w)$
for $\omega\in\Omega,
 w\in W$.

 It is known that (see [Shi])
$c_0=\{w\in\widetilde{W}\mid \textbf{a}(w)=l(w_0)\}$ is a two-sided
cell, which is the lowest one for the partial order $\leq_{LR},$
where $w_0$ is the longest element of $W_0.$ We call $c_0$ the
lowest two-sided cell of $\widetilde{W}.$

In [X1], Xi gave a description of $c_0$ (also can be found in [SX]).
Let $R^+$(resp. $R^-,\Delta$) be the set of positive (resp.
negative, simple) roots in the root system $R$ of $W_0.$ The
dominant weights set $\Lambda^+$ is the set $\{\lambda\in\Lambda\mid
l(\lambda w_0)=l(\lambda)+l(w_0)\}.$ For each simple root $\alpha$
we denote by $s_\alpha$ the corresponding simple reflection in $W_0$
and $x_\alpha$ the corresponding fundamental weight. For each $w\in
W_0,$ we set
$$d_w=w\prod_{{\alpha\in\Delta}\atop {w(\alpha)\in R^-}}x_\alpha.$$
Then $$c_0=\{d_w\lambda w_0d_u^{-1}\mid w,u\in W_0,
\lambda\in\Lambda^+\}.$$

We know that $c_{0,u}=\{d_w\lambda w_0d_u^{-1}\mid w\in W_0,
\lambda\in\Lambda^+\}$ is a left cell in $c_0$ for any $u\in W_0$,
and $c'_{0,w}=\{d_w\lambda w_0d_u^{-1}\mid u\in W_0,
\lambda\in\Lambda^+\}$ is a right cell in $c_0$ for any $w\in W_0$.

Bernstein described the center of Hecke algebras (see [L3]). To each
$\lambda\in \Lambda^+$, Bernstein associates an element
$S_{\lambda}$ in the center of $\widetilde{\mathcal {H}}$. He showed
that $\{S_{\lambda},\lambda\in\Lambda^+\}$ form an $\mathcal
{A}$-basis of the center of $\widetilde{\mathcal {H}}$.

For any $\lambda\in \Lambda^+$, we denote by $V(\lambda)$ a rational
irreducible $G$-module of highest weight $\lambda$.

\noindent {\bf Lemma 2.1} ([ L3, X1, SX]) We have that \\
(1)
$S_{\lambda}S_{\lambda'}=\sum_{z\in\Lambda^+}m_{\lambda,\lambda',z}S_z$
for any $\lambda,\lambda'\in\Lambda^+$. Here
$m_{\lambda,\lambda',z}$ is the multiplicity
of $V(z)$ in the tensor product $V(\lambda)\otimes V(\lambda').$\\
(2) Let $w,w',u\in W_0,\ \lambda,\lambda'\in\Lambda^+$. Then
$$\mu(d_w\lambda
w_0d_u^{-1},d_{w'}\lambda'w_0d_u^{-1})=\mu(d_w\lambda
w_0,d_{w'}\lambda'w_0) .$$ (3) For any elements $u,w \in c_0$ such
that $\widetilde{\mu}(u,w)\neq0$, we have that
$u\sim_Lw$ and $u\nsim_Rw$, or $u\sim_Rw$ and $u\nsim_Lw$. \\
(4) Let $w=d_u\lambda w_0,\ w'=d_{u'}\lambda'w_0,$ where $u,u'\in
W_0,\ \lambda,\lambda'\in\Lambda^+.$ Then the set $\{z\in
\widetilde{W}\mid$ $h_{w_0d_u^{-1},d_{u'}w_0,z}\neq0\}$ is finite,
and if $h_{w_0d_u^{-1},d_{u'}w_0,z}\neq0$ then $z=z_1w_0$ for some
$z_1\in\Lambda^+.$ Moreover
$\mu(w,w')=\sum_{z_1\in\Lambda^+}m_{\lambda^*,\lambda',z_1^{*}}
\delta_{w_0d_u^{-1},d_{u'}w_0,z_1w_0},$ where
$\lambda^*=w_0\lambda^{-1}w_0\in\Lambda^+.$

\section{Cells in an affine Weyl group of type $\widetilde{A_2}$}

In the rest of this paper we assume that
$G=\textrm{SL}_3(\mathbb{C})$. Then we get that $(W,S)$ is an affine
Weyl group of type $\widetilde{A_2}$, where $S=\{r,s,t\}$ satisfies
$(rs)^3=(st)^3=(rt)^3=1$. And $W_0$ is the Weyl group of $G$
generated by $s$ and $t$.

We also know that the corresponding extended affine Weyl group
$\widetilde{W}=\Omega\ltimes W$, where $\Omega$ is the cyclic group
$\{e,\omega,\omega^2\}$ generated by $\omega$ satisfying
$r\omega=\omega s,\ s\omega=\omega t,\ t\omega=\omega r$. We have
that $\omega^{-1}=\omega^2,\ r\omega^2=\omega^2t,\
s\omega^2=\omega^2r$, and $t\omega^2=\omega^2s$.
 The dominant weights set $\Lambda^+=\{\lambda\in\Lambda|\
l(\lambda w_0)=l(\lambda)+l(w_0)\}$, where $w_0=sts$ is the longest
element in $W_0$. Assume that the two simple roots corresponding to
$s,t$ are $\alpha$ and $\beta$, respectively. We have
$$s(\alpha)=-\alpha,\ s(\beta)=t(\alpha)=\alpha+\beta,\
\textrm{and}\ t(\beta)=-\beta.$$ Then the two fundamental dominant
weights of $W$ are
$$x=\frac{2}{3}\alpha+\frac{1}{3}\beta\ \textrm{and}\
y=\frac{1}{3}\alpha+\frac{2}{3}\beta,$$ the corresponding elements
in $\widetilde{W}$ are $$x=\omega tr\ \textrm{and}\ y=\omega^2 sr.$$

We have that $\Lambda=\mathbb{Z}x+\mathbb{Z}y$,
$\Lambda^+=\mathbb{N}x+\mathbb{N}y$, and
$\Lambda_r=\mathbb{Z}\alpha+\mathbb{Z}\beta$.

In [L1], Lusztig described the left cells decomposition of $(W,S)$.
For any subset $J$ of $S=\{r, s, t\}$ we denote by $W^J=\{w\in
W|R(w)=J\}$.  In the following we use the elements in $J$ to denote
$J$. Then $(W,S)$ has $10$ left cells:

$A_{rs}=W^{rs},A_{st}=W^{st},A_{rt}=W^{rt},A_r=A_{st}r,A_s=A_{rt}s,
A_t=A_{rs}t,$

$B_r=W_r-A_r,B_s=W_s-A_s,B_t=W_t-A_t,C_{\emptyset}=W^{\emptyset}=\{e\}.$\\

Set $c_e=C_{\emptyset},c_1=B_t\cup B_s\cup B_t,c_0=A_{rs}\cup
A_{st}\cup A_{rt}\cup A_r\cup A_s\cup A_t$. From [L1], we know that
$c_e,c_1,c_0$ are the entire two-sided cells of $W$. We have that
\begin{eqnarray*}c_e
&=&\{w\in W\mid \textbf{a}(w)=0\}=\{e\}
\\
c_1&=&\{w\in W\mid \textbf{a}(w)=1\}\\
c_0&=&\{w\in W\mid \textbf{a}(w)=3\}.
\end{eqnarray*}

By Lemma 1.1 (7-8), we can easily get that
$${\mu
(e,w)}=\left\{\begin{array}{ll}
{1}, & \ \ \textrm{if} \;  w\in S ,\\
{0} ,   & \ \ \textrm{otherwise}.\end{array} \right.$$ We also have
the following results:

\noindent {\bf Lemma 3.1} ([W1, Proposition 3.4]) For any
$u,w\in c_1,\ u\leq w, $ we have\\
(1) If $L(u)\neq L(w)$ or $R(u)\neq R(w)$, then
$${\mu (u,w)}=\left\{\begin{array}{ll}
{1}, & \ \ \textrm{if} \;  l(w)-l(u)=1,\\
{0} ,   & \ \ \textrm{otherwise};\end{array} \right.$$ (2) If $L(u)=
L(w)$ and $R(u)= R(w)$, then
$${\mu (u,w)}=\left\{\begin{array}{ll}
{1}, & \ \ \textrm{if} \;  l(w)-l(u)=3,\\
{0} ,   & \ \ \textrm{otherwise}.\end{array} \right.$$

The only coefficients $\mu(u,w)$ that we have not known are those
for $u<w$ such that (1) $u,w\in c_0$, (2) $u\in c_1, w\in c_0$, and
(3) $u\in c_0,w\in c_1$. In this paper, we compute $\mu(u,w)$ in
these cases.

The following results will be useful.

\noindent {\bf Lemma 3.2} ([L4, Proposition 3.8(c)]) If $c_1$ is the
second highest two-sided cell for an affine Weyl group, then\\
\noindent (1)
\begin{eqnarray*}
c_1&=&\{w\in W\mid\textbf{a}(w)=1\}\\
&=&\{e\neq w\in W\mid w\ \textrm{has a unique
 reduced\ expression}\}.
\end{eqnarray*}
(2) For $w\in c_1$, we have $|L(w)|=|R(w)|=1$.

\section{The case for $u<w$ and $w\in c_0$}
In this section we compute the Kazhdan-Lusztig coefficients
$\mu(u,w)$ for $u<w$ and $w\in c_0$.

First, we need to describe the lowest two-sided cell of
$\widetilde{W}$ following Section 2. We have
$W_0=\{e,s,t,st,ts,sts=w_0\}$.

\noindent {\bf Proposition 4.1} We have that $d_{e}=e,\ d_{s}=\omega
r,\ d_{t}=\omega^2 r,\ d_{st}=\omega^2 ,\ d_{ts}=\omega,$ and
$d_{sts}=r.$

\noindent {\bf Proof.} From Section 2, we know that
$d_w=w\prod_{{\alpha\in\Delta}\atop {w(\alpha)\in R^-}}x_\alpha$ for
each $w\in c_o$. Obviously, we have that $d_{e}=e$.

Since $s(\a)=-\a$ and $s(\b)=\a+\b$, we get that $d_{s}=s(x)=s\omega
tr=\omega ttr=\omega r$. Also, by facts that $t(\a)=\a+\b$ and
$t(\b)=-\b$ we get that $d_{t}=t(y)=t\omega^2 sr=\omega^2ss
r=\omega^2 r$.

Since $st(\a)=\b$ and $st(\b)=-\a-\b$, we get that
$d_{st}=st(y)=st\omega^2 sr=\omega^2 rssr=\omega^2$.

Since $ts(\a)=-\a-\b$ and $ts(\b)=\a$, we get that
$d_{ts}=ts(x)=ts\omega tr=\omega rttr=\omega$.

Since $sts(\a)=-\a$ and $sts(\b)=-\b$, we get that
$d_{sts}=sts(xy)=sts\omega tr\omega^2 sr=stsstsr=r$. \hfill$\Box$\

Following Lemma 2.1(4), we see that, in order to compute $\mu(w,w')$
for $w<w'$ in $c_0$, we need to know those integers
$\delta_{w_0d_u^{-1},d_{u'}w_0,z_1w_0}\neq0$ for $u,u'\in W_0$ and
$z_1\in\Lambda^+$, where
$\delta_{w_0d_u^{-1},d_{u'}w_0,z_1w_0}\neq0$ is the coefficient of
the power $q$ in $h_{w_0d_u^{-1},d_{u'}w_0,z_1w_0}$. We first
compute those products $C_{w_0d_u^{-1}}C_{d_{u'}w_0}$. Define
$[2]=q^{\frac{1}{2}}+q^{-\frac{1}{2}}.$

\noindent {\bf Proposition 4.2} We have \\
(1) $C_{w_0}C_{w_0}=C_{w_0\omega}C_{\omega^2w_0}=C_{w_0\omega^2}C_{\omega w_0}=([2]^3-[2])C_{w_0}$,\\
(2)
$C_{w_0}C_{rw_0}=C_{w_0r}C_{w_0}=C_{w_0\omega}C_{\omega^2rw_0}=C_{w_0\omega^2}C_{\omega
rw_0}
=C_{w_0r\omega}C_{\omega^2w_0}\\=C_{w_0r\omega^2}C_{\omega w_0}=C_{xyw_0}+[2]^2C_{w_0}$,\\
(3) $C_{w_0}C_{\omega w_0}=C_{w_0\omega}C_{w_0}=C_{w_0\omega^2}C_{\omega^2w_0}=[2]C_{xw_0}$,\\
(4) $C_{w_0}C_{\omega rw_0}=C_{w_0\omega}C_{rw_0}=C_{w_0\omega^2}C_{\omega^2rw_0}=[2]^2C_{xw_0}$,\\
(5) $C_{w_0}C_{\omega^2w_0}=C_{w_0\omega}C_{\omega w_0}=C_{w_0\omega^2}C_{w_0}=[2]C_{yw_0}$,\\
(6) $C_{w_0}C_{\omega^2rw_0}=C_{w_0\omega}C_{\omega
rw_0}=C_{w_0\omega^2}C_{rw_0}=[2]^2C_{yw_0}$.

\noindent {\bf Proof.} By Lemma 1.1 and 1.2, we have that
\begin{eqnarray*}
C_{w_0}&=&C_{sts} =C_{s}C_{ts}-C_{s} =C_{s}(C_{t}C_{s}-1)
=(C_{s}C_{t}-1)C_{s}\\
&=&C_{tst} =C_{t}C_{st}-C_{t} =C_{t}(C_{s}C_{t}-1)
=(C_{t}C_{s}-1)C_{t}.
\end{eqnarray*}
(1) We have
\begin{eqnarray*}
C_{w_0}C_{w_0}&=&C_{sts}C_{sts}
=(C_{s}C_{t}-1)C_{s}C_{sts}\\
&=&[2](C_{s}C_{t}C_{sts}-C_{sts})\\
&=&([2]^3-[2])C_{w_0}.
\end{eqnarray*}
Since $\omega^3=e$, we get that
$C_{w_0\omega}C_{\omega^2w_0}=C_{w_0}C_{\omega^3}C_{w_0}=C_{w_0}C_{w_0}$
and\\ $C_{w_0\omega^2}C_{\omega w_0}
=C_{w_0}C_{\omega^3}C_{w_0}=C_{w_0}C_{w_0}$.

\noindent(2) We have
\begin{eqnarray*}
C_{w_0}C_{rw_0}&=&C_{sts}C_{rsts}\\
&=&(C_{s}C_{t}-1)C_{s}C_{rsts}\\
&=&(C_{s}C_{t}-1)(C_{srsts}+C_{sts})\\
&=&C_{s}C_{tsrsts}+[2]^2C_{w_0}-C_{srsts}-C_{w_0}\\
&=&C_{stsrsts}+C_{srsts}+C_{w_0}+[2]^2C_{w_0}-C_{srsts}-C_{w_0}\\
&=&C_{xyw_0}+[2]^2C_{w_0}.
\end{eqnarray*}
Obviously, we have
$C_{w_0r\omega}C_{\omega^2w_0}=C_{w_0r\omega^2}C_{\omega
w_0}=C_{w_0r}C_{w_0}=C_{w_0}C_rC_{w_0}=C_{w_0}C_{rw_0}=C_{w_0\omega}C_{\omega^2rw_0}=C_{w_0\omega^2}C_{\omega
rw_0}$.

\noindent(3) We have
\begin{eqnarray*}
C_{w_0}C_{\omega w_0}&=&(C_{s}C_{t}-1)C_{s}C_{rsr\omega}\\
&=&[2](C_{s}C_{t}-1)C_{rsr\omega}\\
&=&[2]C_{s}C_{trsr\omega}-[2]C_{rsr\omega}\\
&=&[2](C_{strsr\omega}+C_{rsr\omega})-[2]C_{rsr\omega}\\
&=&[2]C_{xw_0}.
\end{eqnarray*}
Obviously, we have
$C_{w_0\omega}C_{w_0}=C_{w_0\omega^2}C_{\omega^2w_0}=C_{w_0}C_{\omega
w_0}$.

\noindent(4) We have
\begin{eqnarray*}
C_{w_0}C_{\omega rw_0}&=&(C_{t}C_{s}-1)C_{t}C_{trsr\omega}\\
&=&[2](C_{t}C_{s}-1)C_{trsr\omega}\\
&=&[2]C_{t}(C_{strsr\omega}+C_{rsr\omega})-[2]C_{trsr\omega}\\
&=&[2]^2C_{strsr\omega}+[2]C_{trsr\omega}-[2]C_{trsr\omega}\\
&=&[2]C_{xw_0}.
\end{eqnarray*}
Obviously, we have
$C_{w_0\omega}C_{rw_0}=C_{w_0\omega^2}C_{\omega^2rw_0}=C_{w_0}C_{\omega
rw_0}$.

\noindent(5) We have
\begin{eqnarray*}
C_{w_0}C_{\omega^2 w_0}&=&(C_{t}C_{s}-1)C_{t}C_{trt\omega^2}\\
&=&[2](C_{t}C_{s}-1)C_{trt\omega^2}\\
&=&[2]C_{t}C_{strt\omega^2}-[2]C_{trt\omega^2}\\
&=&[2](C_{tstrt\omega^2}+C_{trt\omega^2})-[2]C_{trt\omega^2}\\
&=&[2]C_{yw_0}.
\end{eqnarray*}
Obviously, we have  $C_{w_0\omega}C_{\omega
w_0}=C_{w_0\omega^2}C_{w_0}=C_{w_0}C_{\omega^2w_0}$.

\noindent(6) We have
\begin{eqnarray*}
C_{w_0}C_{\omega^2 rw_0}&=&(C_{s}C_{t}-1)C_{s}C_{strt\omega^2}\\
&=&[2](C_{s}C_{t}-1)C_{strt\omega^2}\\
&=&[2]C_{s}(C_{tstrt\omega^2}+C_{trt\omega^2})-[2]C_{strt\omega^2}\\
&=&[2]^2C_{tstrt\omega^2}+[2]C_{strt\omega^2}-[2]C_{strt\omega^2}\\
&=&[2]^2C_{yw_0}.
\end{eqnarray*}
Obviously, we have  $C_{w_0\omega}C_{\omega
rw_0}=C_{w_0\omega^2}C_{rw_0}=C_{w_0}C_{\omega^2rw_0}$.
\hfill$\Box$\

\noindent {\bf Proposition 4.3} We have \\
(1)
$C_{w_0r}C_{rw_0}=C_{w_0r\omega}C_{\omega^2rw_0}=C_{w_0r\omega^2}C_{\omega
rw_0}=
[2]C_{xyw_0}+[2]^3C_{w_0}$,\\
(2) $C_{w_0r}C_{\omega w_0}=C_{w_0r\omega}C_{w_0}=C_{w_0r\omega^2}C_{\omega^2w_0}=[2]^2C_{xw_0}$,\\
(3) $C_{w_0r}C_{\omega^2w_0}=C_{w_0r\omega^2}C_{w_0}=C_{w_0r\omega}C_{\omega w_0}=[2]^2C_{yw_0}$,\\
(4) $C_{w_0r}C_{\omega
rw_0}=C_{w_0r\omega}C_{rw_0}=C_{w_0r\omega^2}C_{\omega^2rw_0}=[2]C_{y^2w_0}+2[2]C_{xw_0}$,\\
(5)
$C_{w_0r}C_{\omega^2rw_0}=C_{w_0r\omega^2}C_{rw_0}=C_{w_0r\omega}C_{\omega
rw_0}=[2]C_{x^2w_0}+2[2]C_{yw_0}$.

\noindent {\bf Proof.} The proofs are very similar to those of
Proposition 4.2.

\noindent(1) We have
\begin{eqnarray*}
C_{w_0r}C_{rw_0}&=&C_{w_0}C_{r}C_{rsts}\\
&=&[2]C_{w_0}C_{rsts}\\
&=&[2]C_{xyw_0}+[2]^3C_{w_0}.
\end{eqnarray*}
Obviously, we have
$C_{w_0r\omega}C_{\omega^2rw_0}=C_{w_0r\omega^2}C_{\omega
rw_0}=C_{w_0r}C_{ rw_0}$.

\noindent(2) We have
\begin{eqnarray*}
C_{w_0r}C_{\omega w_0}&=&C_{w_0}C_{r}C_{rsr\omega}\\
&=&[2]C_{w_0}C_{\omega w_0}\\
&=&[2]^2C_{xw_0}.
\end{eqnarray*}
Obviously, we have
$C_{w_0r\omega}C_{w_0}=C_{w_0r\omega^2}C_{\omega^2w_0}=C_{w_0r}C_{\omega
w_0}$.

\noindent(3) We have
\begin{eqnarray*}
C_{w_0r}C_{\omega^2 w_0}&=&C_{w_0}C_{r}C_{trt\omega^2}\\
&=&[2]C_{w_0}C_{\omega^2 w_0}\\
&=&[2]^2C_{yw_0}.
\end{eqnarray*}
Obviously, we have $C_{w_0r\omega^2}C_{w_0}=C_{w_0r\omega}C_{\omega
w_0}=C_{w_0r}C_{\omega^2 w_0}$.

\noindent(4) We have
\begin{eqnarray*}
C_{w_0r}C_{\omega rw_0}&=&C_{w_0}C_{r}C_{trsr\omega}\\
&=&(C_{t}C_{s}-1)C_{t}(C_{rtrsr\omega}+C_{rsr\omega})\\
&=&[2](C_{t}C_{s}-1)C_{rtrsr\omega}+C_{w_0}C_{\omega w_0}\\
&=&[2]C_{t}C_{srtrsr\omega}-[2]C_{rtrsr\omega}+C_{w_0}C_{\omega w_0}\\
&=&[2](C_{tsrtrsr\omega}+C_{rtrsr\omega}+C_{strsr\omega})-
[2]C_{rtrsr\omega}+[2]C_{xw_0}\\
&=&[2]C_{y^2w_0}+2[2]C_{xw_0}.
\end{eqnarray*}
Obviously, we have
$C_{w_0r\omega}C_{rw_0}=C_{w_0r\omega^2}C_{\omega^2rw_0}=C_{w_0r}C_{\omega
rw_0}$.

\noindent(5) We have
\begin{eqnarray*}
C_{w_0r}C_{\omega^2 rw_0}&=&C_{w_0}C_{r}C_{strt\omega^2}\\
&=&C_{w_0}(C_{rstrt\omega^2}+C_{trt\omega^2})\\
&=&(C_{s}C_{t}-1)C_{s}C_{rstrt\omega^2}+C_{w_0}C_{\omega^2 w_0}\\
&=&[2](C_{s}C_{t}-1)C_{rstrt\omega^2}+C_{w_0}C_{\omega^2 w_0}\\
&=&[2]C_{s}C_{trstrt\omega^2}-[2]C_{rstrt\omega^2}+C_{w_0}C_{\omega^2 w_0}\\
&=&[2](C_{strstrt\omega^2}+C_{rstrt\omega^2}+C_{tstrt\omega^2})-[2]C_{rstrt\omega^2}+[2]C_{yw_0} \\
&=&[2]C_{x^2w_0}+2[2]C_{yw_0}.
\end{eqnarray*}
Obviously, we have $C_{w_0r\omega^2}C_{rw_0}=C_{w_0r\omega}C_{\omega
rw_0}=C_{w_0r}C_{\omega^2 rw_0}$. \hfill$\Box$\

\noindent {\bf Proposition 4.4} For any $u,u'\in W_0,$ there is at
most one $z_1\in\Lambda^+$ such that
$\delta_{w_0d_u^{-1},d_{u'}w_0,z_1w_0}\neq0.$ Moreover, if
$\delta_{w_0d_u^{-1},d_{u'}w_0,z_1w_0}\neq0,$ we get
$\delta_{w_0d_u^{-1},d_{u'}w_0,z_1w_0}=1$ and the unique $z_1$ is in
$\{0,x,y\},$ where $x,y\in\Lambda^+$ are the fundamental dominant
weights.

\noindent {\bf Proof.} By the definition of
$\delta_{w_0d_u^{-1},d_{u'}w_0,z_1w_0}$ in Section 1 and
$\textbf{a}(z_1w_0)=3$, we get that
$\delta_{w_0d_u^{-1},d_{u'}w_0,z_1w_0}$ is the coefficient of the
power $q$ in $h_{w_0d_u^{-1},d_{u'}w_0,z_1w_0}$. Combing
Propositions 4.2 and 4.3, we can get the results. \hfill$\Box$\

If we define $\mathcal {U}$ $= \{(d_u,d_{u'},z)\mid u,u'\in
W_0,z\in\Lambda^{+},\delta_{w_0d_u^{-1},d_{u'}w_0,zw_0}=1\}$, then
by Propositions 4.2 and 4.3 we get that $\mathcal
{U}=\{(e,r,0),(e,\omega r,x),(e,\omega^2r,y),\\(r,e,0),(r,\omega,x),
(r,\omega^2,y),(\omega,r,x),(\omega,\omega
r,y),(\omega,\omega^2r,0), (\omega^2,r,y),(\omega^2,\omega
r,0),\\(\omega^2,\omega^2r,x), (\omega r,e,x),(\omega
r,\omega,y),(\omega r,\omega^2,0),(\omega^2r,e,y),
(\omega^2r,\omega,0),(\omega^2r,\omega^2,x),\}$

Following Lemma is a well-known result in representation theory (see
[BZ]):

\noindent {\bf Lemma 4.5} Assume that $L$ is a semisimple complex
Lie algebra, $\lambda\ \textrm{and}\ \lambda'$ are dominant weights.
Let $K_{\lambda,\beta}$ be the weight multiplicity of weight $\beta$
in the irreducible $L$-module $V(\lambda),$ and
$m_{\lambda,\lambda',\nu}$ be the tensor product multiplicity of
irreducible module $V(\nu)$ in $V(\lambda)\otimes V(\lambda').$ Then
we have that $0\leq m_{\lambda,\lambda',\nu}\leq
K_{\lambda,\nu-\lambda'}.$

Now we get the main results in this section:

\noindent {\bf Theorem 4.6} For any $w,w'\in c_0$, if $w\sim_Lw'$
i.e. $w=d_u\lambda w_0d_v^{-1},\ w'=d_{u'}\lambda'w_0d_v^{-1},$ for
some $u,u',v\in W_0$ and $\lambda,\lambda'\in \Lambda^{+},$ then we
have that $\mu(w,w')=m_{\lambda^*,\lambda',z_1^*}$ if there is some
$z_1\in\{0,x,y\}$ such that $(d_u,d_{u'},z_1)\in$ $\mathcal {U}$;
otherwise, $\mu(w,w')=0.$ Moreover, we get that for any $w,w'\in
c_0,$ $\mu(w,w')\leq1.$

\noindent {\bf Proof.} For any $w,w'\in c_0,$ if $\mu(w,w')\neq0,$
we have $w\sim_Lw'$ or $w\sim_Rw'$ by Lemma 1.3. By Lemma 1.1 (6)
and Lemma 1.4, we can assume that $w\sim_Lw'.$ If $w=d_uxw_0d_v,\
w'=d_{u'}x'w_0d_v,$ for some $u,u',v\in W_0$ and
$\lambda,\lambda'\in \Lambda^{+},$ by Lemma 2.1 (2), we can also
assume that $w=d_uxw_0,\ w'=d_{u'}x'w_0.$ From Lemma 2.1 (4), we
know that
$$\mu(w,w')=\sum_{z_1\in\Lambda^+}m_{\lambda^*,\lambda',z_1^{*}}
\delta_{w_0d_u^{-1},d_{u'}w_0,z_1w_0},\ \textrm{where}\
\lambda^*=w_0\lambda^{-1}w_0\in\Lambda^+.$$

We can check that $x^*=y$ and $y^*=x$ in the case of type $A_2$.

By Proposition 4.4, we get that
$\mu(w,w')=m_{\lambda^*,\lambda',z_1^{*}}$ if there is some
$z_1\in\{0,x,y\}$ such that $(d_u,d_{u'},z_1)\in$ $\mathcal {U}$;
otherwise, $\mu(w,w')=0.$

On the other side, we have\begin{eqnarray*}(*)\ \ \ \ \ \
m_{\lambda^*,\lambda',z_1^{*}}&=&\textrm{dim\
Hom}_\mathcal {G}(V(\lambda^*)\otimes V(\lambda'),V(z_1^*))\\
&=&\textrm{dim\ Hom}_\mathcal {G}(V(z_1)\otimes
V(\lambda'),V(\lambda))\\
&=&m_{z_1,\lambda',\lambda},\end{eqnarray*}where $\mathcal {G}$ is
the Lie algebra of $G$.

Since $x=\frac{1}{3}(2\alpha+\beta)$ and
$y=\frac{1}{3}(\alpha+2\beta)$, thus the weight set of $V(z_1)$
equals to $W_0(z_1)$ if $z_1\in\{0,x,y\}$. Then
$K_{z_1,\lambda}=K_{z_1,z_1}=1$ for any weight $\lambda$ of
$V(z_1).$ By Lemma 4.5, we get
$m_{\lambda^*,\lambda',z_1^{*}}=m_{z_1,\lambda',\lambda}\leq
K_{z_1,\lambda-\lambda'}\leq1$ for $\lambda,\lambda'\in\Lambda^+$.
The theorem holds. \hfill$\Box$\

\noindent {\bf Corollary 4.7} For any elements $w\leq w'\in c_0$, if
$\mu (w,w')\neq0$, then we have that $l(w')-l(w)=1\ \textrm{or}\ 3$.

\noindent {\bf Proof.} Assume that $z_1=0,\ x$ or $y$. We can assume
that $w=d_u\lambda w_0d_v<\ w'=d_{u'}\lambda'w_0d_v$ for some
$u,u',v\in W_0$ and $\lambda,\lambda'\in \Lambda^+$. By Theorem 4.6
and the equalities ($*$) in its proof, we know that
$\mu(w,w')=m_{z_1,\lambda',\lambda}$. By the theory of Lie algebra,
we know that  $m_{z_1,\lambda',\lambda}\neq0$ implies that
$\lambda-\lambda'$ must be a weight  of $V(z_1)$.

If $z_1=0$, then $m_{0,\lambda',\lambda}\neq0$ implies that
$\lambda=\lambda'$.

If $z_1=x$, then $m_{x,\lambda',\lambda}\neq0$ implies that
$\lambda-\lambda'\in\{x,y-x,-y\}$. With the assumption $w<w'$, we
get that $\lambda-\lambda'=-y$.

If $z_1=y$, then $m_{y,\lambda',\lambda}\neq0$ implies that
$\lambda-\lambda'\in\{y,-x,x-y\}$. With the assumption $w<w'$, we
get that $\lambda-\lambda'=-x$.

Since $l(x)=l(y)=2$,  we see that $\mu(w,w')\neq0$ implies that
$l(\lambda')-l(\lambda)=$0 or 2. Then we get the result, since
$l(d_u)=$0 or 1 for any $u\in W_0$. \hfill$\Box$\

Now we compute those $\mu(u,w)$ for $u<w$ such that $u\in c_1$ and
$w\in c_0$. By [X2], we know that $\mu(u,w)\leq1$ in this case. Here
we give these values explicitly.

\noindent {\bf Theorem 4.8} For elements $u<w\in W$ such that $u\in
c_1$ and $w\in c_0$, we have that
$${\mu(u,w)}=\left\{\begin{array}{ll}
{1},& \ \ \textrm{if} \ l(w)-l(u)=1,\\
{0},& \ \ \textrm{otherwise}.\end{array} \right.$$

\noindent {\bf Proof.} Assume that $\mu(u,w)\neq0$. Then by Lemma
1.3, we get that $w\leq_Ru$ and $w\leq_Lu$. Thus $L(u)\subseteq
L(w)$ and $R(u)\subseteq R(w)$ by Lemma 1.4.

If $L(u)\subsetneq L(w)$ or $R(u)\subsetneq R(w)$, we get the result
by Lemma 1.1 (7--8).

Now we assume that $L(u)= L(w)$ and $R(u)= R(w)$. Without loss
generality, we can assume that $L(u)=\{r\}$. This follows that
$L(w)=\{r\}$ and $L(rw)=\{s,t\}$. Assume that $L(ru)=\{s_1\}$ for
some $s_1\in\{s,t\}$. We denote by $s_2\in\{s,t\}\backslash
\{s_1\}$. Then we use star operation, which is defined in [KL1], on
the left with respect to $\{r,s_2\}$ and by [KL1, Theorem 5.2]. We
get that $\mu(u,w)=\mu(s_2u,rw)$. This implies that
$l(w)-l(u)=l(rw)-l(s_2u)-2$ and $L(s_2u)=\{s_2\}$, $L(rw)=\{s,t\}$.
If $l(w)-l(u)\geq3$, then by Lemma 1.1 (8), we know that
$\mu(s_2u,rw)\neq0$ if and only if $s_1s_2u=rw$. But this cannot
hold, because $s_1s_2u\in c_1$ and $rw\in c_0$. Thus we get that
$\mu(u,w)=0$, if $l(w)-l(u)\geq3$ in the case that $L(u)= L(w)$ and
$R(u)= R(w)$. If $l(w)-l(u)=1$, by Lemma 1.1 (2), we get that
$\mu(u,w)=1$.

We complete the proof. \hfill$\Box$\

In conclusion, we see that if $u\leq w$ and $\textbf{a}(u)\leq
\textbf{a}(w)$, then $\mu(u,w)\leq1$ and if $\mu(u,w)\neq0$, then
$l(w)-l(u)=1$ or $3$.

\section{The case for $(u,w)\in c_0\times c_1$ and $l(w)-l(u)\leq3$}
For any $w\in W$, $E\subseteq W$, and $n\in \mathbb{N}$, we define a
subset $E_n^w=\{u\in E\mid u<w\ \textrm{and}\ l(w)-l(u)=n\}$.

We consider $E=c_0$. First, we
 describe the sets $(c_0)_1^w$ for $w\in c_1$. By Lemma 3.2 (1), the
 reduced expression of $w$ is unique. Obviously, any element in
 $(c_0)_1^w$ can be obtained from the reduced
expression of $w$ by removing one proper generator. Recall our
Notations (2) in section Introduction.

\noindent {\bf Lemma 5.1} If the reduced expression of $w\in c_1$ is
$s_1s_2\cdots s_{n}$, then any $u\in (c_0)_1^w$ equals to
$s_{1}\widehat{s}_{2}s_{3}\cdots s_{n}$ or $s_{1}s_{2}\cdots
s_{n-2}\widehat{s}_{n-1}s_{n}$.

\noindent {\bf Lemma 5.2} Let $w\in c_1$ and $u\in (c_0)_1^w$. If
$R(w)=\{s'\}$ and $|L(u)|=1$, then $R(ws')\nsubseteq R(u)$.

\noindent {\bf Proof.} Without loss of generality, we assume that
$w=(rst)^mw'$, where $m\in \mathbb{N}$ and $w'\in\{e,r,rs\}$.

(1) If $w=(rst)^m$, by Lemma 5.1 and $|L(u)|=1$, we get that
$u=(rst)^{m-1}rt=rs(trs)^{m-2}trt$ and $m$ is odd. Since
$R(w)=\{t\}$, $R(wt)=\{s\}$, and $R(u)=\{r,t\}$, we have
$R(wt)\nsubseteq R(u)$.

(2) If $w=(rst)^mr$, by Lemma 5.1 and $|L(u)|=1$, we get that
$u=(rst)^{m-1}rsr$ and $m$ is even. Since $R(w)=\{r\}$,
$R(wr)=\{t\}$, and $R(u)=\{r,s\}$, we have $R(wr)\nsubseteq R(u)$.

(3) If $w=(rst)^mrs$, by Lemma 5.1 and $|L(u)|=1$, we get that
$u=(rst)^{m-1}rsts$ and $m$ is odd. Since $R(w)=\{s\}$,
$R(wr)=\{r\}$, and $R(u)=\{s,t\}$, we have $R(ws)\nsubseteq R(u)$.

The Lemma is proved.  \hfill$\Box$\

Similarly, we have

\noindent {\bf Lemma 5.3} Let $w\in c_1$ and $u\in (c_0)_1^w$. If
$L(w)=\{s'\}$ and $|R(u)|=1$, then $L(s'w)\nsubseteq L(u)$.

In the following we fix an element $w\in c_1$, whose unique reduced
expression is denoted by $w=s_1s_2\cdots s_{n}$ satisfying
$\{s_1,s_2,s_3\}=S=\{r,s,t\}$. By Lemma 3.2 (2), we know that
$L(w)=\{s_1\}$ and $R(w)=\{s_n\}$. For some $u\in (c_0)_3^w$, we
define

\noindent {\bf Condition 5.4} $|L(u)|=|R(u)|=2$, $s_1\in L(u)$,
$s_n\in R(u)$, $\{s_2\}\nsubseteq L(u)$, and $\{s_{n-1}\}\nsubseteq
R(u)$.

First we have two propositions.

\noindent {\bf Proposition 5.5} Let $w\in c_1$ and $u\in (c_0)_3^w$.
If $|L(u)|=|R(u)|=2$, $s_1\in L(u)$, $s_n\in R(u)$, $s_2\in L(u)$,
or $s_{n-1}\in R(u)$, then $\mu(u,w)=0$.

\noindent {\bf Proof.} Element $u\in (c_0)_3^w$ can be gotten from
the reduced expression $w=s_1s_2\cdots s_{n}$ by deleting three
proper letters. In order to satisfy $|L(u)|=|R(u)|=2$, $s_1\in L(u)$
and $s_n\in R(u)$, conditions $s_2\in L(u)$ and $s_{n-1}\in R(u)$
can not hold at the same time. We can assume that $s_{n-1}\in R(u)$,
but $s_2$ is not in $L(u)$.

Hence we must have $w=s_1s_2s_3w'$ and $u=s_1s_3s_1u'$ such that
$l(w)=l(w')+3$ and $l(u)=l(u')+3$. Then we use star operation, which
is defined in [KL1], on the left with respect to $\{s_1,s_2\}$ and
by [KL1, Theorem 5.2]. We get that $\mu(u,w)=\mu(s_2u,s_1w)$. Also,
we see that $s_1w\in c_1$, $l(s_1w)-l(s_2u)=[l(w)-1]-[l(u)+1]=1$,
$L(s_2u)=\{s_2\}$, $R(s_1w)=R(w)$ and $R(s_2u)=R(u)$; that is,
$R(s_1ws_n)\subseteq R(s_2u)$. By Lemma 5.2, we see $s_2u$ is not in
$(c_0)1^{s_1w}$, which implies that $s_2u\nleq s_1w$. Thus, we get
that $\mu(u,w)=\mu(s_2u,s_1w)=0$. \hfill$\Box$\

\noindent {\bf Proposition 5.6} Let $w\in c_1$ and $u\in (c_0)_3^w$.
If $|L(u)|=1$ or $|R(u)|=1$, $s_1\in L(u)$ and $s_n\in R(u)$, then
$\mu(u,w)=0$.

\noindent {\bf Proof.} Without loss of generality, we assume that
$L(u)=\{s_1\}$. We must have that $u=s_1s_2s_3s_2u'$ and
$w=s_1s_2s_3w'$ such that $l(u)=l(u')+4$ and $l(w)=l(w')+3$. Using
star operation on the left with respect to $\{s_1,s_2\}$ ([KL1,
Theorem 5.2]), we get $\mu(u,w)=\mu(s_1u,s_1w)$. We see that
$s_1u\in c_0$, $s_1w\in c_1$, $l(s_1w)-l(s_1u)=[l(w)-1]-[l(u)-1]=3$,
$|L(s_1u)|=2$, $L(s_1w)=\{s_2\}\subseteq L(s_1u)$, and $s_3\in
L(s_1u)$. Moreover, we have $R(u)=R(s_1u)$.

If $|R(u)|=1$, we can using star operation on the right in the
similar way. Thus we can get $\mu(u,w)=\mu(u'',w'')$, the $u''$ and
$w''$ satisfy the conditions in Proposition 5.5. The proof is
finished.\hfill$\Box$\

Now we get one of the main results in this section.

\noindent {\bf Theorem 5.7} For $w\in c_1$ and $u\in (c_0)_3^w$, if
$\mu(u,w)\neq0$, then  $u$ must satisfy Condition 5.4.

\noindent {\bf Proof.} By Lemma 1.3 and $\mu(u,w)\neq0$, we get
$u\leq_L w$ and $u\leq_R w$. By Lemma 1.4, we have $R(w)\subseteq
R(u)$ and $L(w)\subseteq L(u)$; that is, $s_1\in L(u)$ and $s_n\in
R(u)$. Combine Propositions 5.5 and 5.6, we see that if $u$ does not
satisfy the Condition 5.4, we must have $\mu(u,w)=0$. Thus the
theorem holds. \hfill$\Box$\

In the following, we compute $\mu(u,w)$ for $w\in c_1$ and those
$u\in(c_0)_3^w$ satisfying Condition 5.4. Also, we assume  that
$\{s_1,s_2,s_3\}=S=\{r,s,t\}$.

\noindent {\bf Proposition 5.8} Let $w\in c_1$ is the given element
with the unique reduced expression $w=s_1s_2\cdots s_{n}$. Then
there exists at most one element $u\in(c_0)_3^w$ satisfying
Condition 5.4, and the unique element is in the form
$u=\widehat{s}_1\widehat{s}_2s_3\cdots
\widehat{s}_{n-1}s_{n}=s_1\widehat{s}_2s_3\cdots
s_{n-2}\widehat{s}_{n-1}\widehat{s}_{n}$.

\noindent {\bf Proof.} We can get an element $u\in(c_0)_3^w$ by
deleting three proper generators from the unique reduced expression
$w=s_1s_2\cdots s_{n}$. If $u$ satisfies Condition 5.4 and
$l(u)=l(w)-3$, we get that $u$ must be in the form
$\widehat{s}_1\widehat{s}_2s_3\cdots \widehat{s}_{n-1}s_{n}$ or
$s_1\widehat{s}_2s_3\cdots s_{n-2}\widehat{s}_{n-1}\widehat{s}_{n}$.
In fact, if $u$ satisfies Condition 5.4, we must get
$\widehat{s}_1\widehat{s}_2s_3\cdots
\widehat{s}_{n-1}s_{n}=s_1\widehat{s}_2s_3\cdots
s_{n-2}\widehat{s}_{n-1}\widehat{s}_{n}$.\hfill$\Box$\

\noindent {\bf Proposition 5.9} Let $w=(s_1s_2s_3)^m$, for some
$m\in\mathbb{N}$, and $u\in(c_0)_3^w$. Then $\mu(u,w)\neq0$ if and
only if $u=(s_3s_1s_2)^{m-2}s_3s_1s_3$ and $m\geq2$ is even.
Moreover, if $\mu(u,w)\neq0$, we have $\mu(u,w)=1$.

\noindent {\bf Proof.} By Theorem 5.7, we see $\mu(u,w)\neq0$
implies that $u$ must satisfy Condition 5.4. By Proposition 5.8, we
get that $u=(s_3s_1s_2)^{m-2}s_3s_1s_3=s_3s_1s_3(s_2s_3s_1)^{m-2}$
and $m\geq2$ is even.

Now we assume that $m$ is even. Using left star operation with
respect to $\{s_1,s_2\}$, we get that
\begin{eqnarray*}
\mu(u,w)&=&\mu((s_3s_1s_2)^{m-2}s_3s_1s_3,(s_1s_2s_3)^m)\\
&=&\mu(s_2(s_3s_1s_2)^{m-2}s_3s_1s_3,s_2s_3(s_1s_2s_3)^{m-1})\\
&=&\mu((s_2s_3s_1)^{m-1}s_3,(s_2s_3s_1)^{m-1}s_2s_3)\\
&=&1.
\end{eqnarray*}\hfill$\Box$\

\noindent {\bf Proposition 5.10} Let $w=(s_1s_2s_3)^ms_1$, for some
$m\in\mathbb{N}$, and $u\in(c_0)_3^w$. Then $\mu(u,w)\neq0$ if and
only if $u=(s_3s_1s_2)^{m-2}s_3s_1s_2s_1$ and $m\geq2$ is odd.
Moreover, if $\mu(u,w)\neq0$, we have $\mu(u,w)=1$.

\noindent {\bf Proof.} By Theorem 5.7, we see $\mu(u,w)\neq0$
implies that $u$ must satisfy Condition 5.4. By Proposition 5.8, we
get that
$u=(s_3s_1s_2)^{m-2}s_3s_1s_2s_1=s_3s_1s_3(s_2s_3s_1)^{m-2}s_2$ and
$m\geq2$ is odd.

Now we assume that $m$ is odd. Using left star operation with
respect to $\{s_1,s_2\}$, we get that
\begin{eqnarray*}
\mu(u,w)&=&\mu((s_3s_1s_2)^{m-2}s_3s_1s_2s_1,(s_1s_2s_3)^ms_1)\\
&=&\mu(s_2(s_3s_1s_2)^{m-2}s_3s_1s_2s_1,s_2s_3(s_1s_2s_3)^{m-1}s_1)\\
&=&\mu((s_2s_3s_1)^{m-1}s_2s_1,(s_2s_3s_1)^{m-1}s_2s_3s_1)\\
&=&1.
\end{eqnarray*}\hfill$\Box$\

\noindent {\bf Proposition 5.11} Let $w=(s_1s_2s_3)^ms_1s_2$, for
some $m\in\mathbb{N}$, and $u\in(c_0)_3^w$. Then $\mu(u,w)\neq0$ if
and only if $u=(s_3s_1s_2)^{m-2}s_3s_1s_2s_3s_2$ and $m\geq2$ is
even. Moreover, if $\mu(u,w)\neq0$, we have $\mu(u,w)=1$.

\noindent {\bf Proof.} By Theorem 5.7, we see $\mu(u,w)\neq0$
implies that $u$ must satisfy Condition 5.4. By Proposition 5.8, we
get that
$u=(s_3s_1s_2)^{m-2}s_3s_1s_2s_3s_2=s_3s_1s_3(s_2s_3s_1)^{m-2}s_2s_3$
and $m\geq2$ is even.

Now we assume that $m$ is odd. Using left star operation with
respect to $\{s_1,s_2\}$, we get that
\begin{eqnarray*}
\mu(u,w)&=&\mu((s_3s_1s_2)^{m-2}s_3s_1s_2s_3s_2,(s_1s_2s_3)^ms_1s_2)\\
&=&\mu(s_2(s_3s_1s_2)^{m-2}s_3s_1s_2s_3s_2,s_2s_3(s_1s_2s_3)^{m-1}s_1s_2)\\
&=&\mu((s_2s_3s_1)^{m-1}s_2s_3s_2,(s_2s_3s_1)^{m-1}s_2s_3s_1s_2)\\
&=&1.
\end{eqnarray*}\hfill$\Box$\

By the above results, we can get the following result easily.

\noindent {\bf Theorem 5.12} For $w\in c_1$ and $u\in (c_0)_3^w$, we
have $\mu(u,w)\neq0$ if and only if $u$ satisfies Condition 5.4.
Moreover, the nonzero value is $1$.

\section{The case for $(u,w)\in c_0\times c_1$ and $l(w)-l(u)\geq5$}
In [L3] Lusztig introduced some semilinear equations, which are
useful for calculating some $\mu(u,w)$.
 In this section we first recall the equations, then discuss their relevance
   to  type $\widetilde{A}_2$.

Recall that $(W,S)$ is an affine Weyl group of type
$\widetilde{A_2}$, where $S=\{r,s,t\}$ such that
$(rs)^3=(st)^3=(rt)^3=1$, and $W_0$ is the Weyl group of $G$
generated by $s$ and $t$.

We have $\widetilde{W}=\Omega\ltimes W$, where $\Omega$ is the
cyclic group $\{e,\omega,\omega^2\}$ generated by $\omega$ and
satisfying $r\omega=\omega s,\ s\omega=\omega t,\ t\omega=\omega r$.
We get that $\omega^{-1}=\omega^2,\ r\omega^2=\omega^2t,\
s\omega^2=\omega^2r,\ t\omega^2=\omega^2s$.

 The dominant weights set $\Lambda^+=\{\lambda\in\Lambda|\
l(\lambda w_0)=l(\lambda)+l(w_0)\}$, where $w_0=sts$ is the longest
element in $W_0$. Assume that the two simple roots corresponding to
$s,t$ are $\alpha$ and $\beta$, respectively. We have
$$s(\alpha)=-\alpha,\ s(\beta)=t(\alpha)=\alpha+\beta,\
\textrm{and}\ t(\beta)=-\beta.$$ Then the two fundamental dominant
weights of $W$ are
$$x=\frac{2}{3}\alpha+\frac{1}{3}\beta\ \textrm{and}\
y=\frac{1}{3}\alpha+\frac{2}{3}\beta,$$ the corresponding elements
in $\widetilde{W}$ are $$x=\omega tr\ \textrm{and}\ y=\omega^2 sr.$$

We get that $$\a=2x-y\ \textrm{and}\ \b=-x+2y.$$We have that
$\Lambda=\mathbb{Z}x+\mathbb{Z}y$,
$\Lambda^+=\mathbb{N}x+\mathbb{N}y$, and
$\Lambda_r=\mathbb{Z}\alpha+\mathbb{Z}\beta\subset\Lambda$.

\def\L{\Lambda}
\def\l{\lambda}

 We need some notations. Recall that $R^+$ is the subset of $R$ containing
  all positive roots. For any subset
$\textbf{i}\subseteq R^+,$ we set
$\alpha_{\textbf{i}}=\sum_{\alpha\in \textbf{i}}\alpha \in \Lambda$.
Set $\rho=\frac 12\alpha_{R^+}.$ Then $\rho$ is a dominant weight
and $w(\rho)-\rho\in \Lambda_r$ for any $w\in W_0.$ For any
$\lambda\in \Lambda_r$ we set
$$\Phi(\lambda)=\sum_{{}\atop{\textbf{i}\subseteq R^+};\
\alpha_{\textbf{i}} =\lambda}(-v^2)^{-|\textbf{i}|}.$$ Note that the
summation index $\textbf{i}$ runs through all subsets $\textbf{i}$
of $R^+$ such that $\lambda$ can be written as the sum of all
elements in $\textbf{i}$. If $0\neq\lambda$ cannot be written as a
sum of distinct positive roots,  we set $\Phi(\lambda)=0.$

 For any $\lambda\in
\Lambda^+,$ we set $$W_0^{\lambda}=\{w\in W_0\mid
w(\lambda)=\lambda\}$$ and
$$\pi_{\lambda}=v^{-\nu_{\lambda}}\sum_{w\in
W_0^{\lambda}}v^{2l(w)},$$ where $\nu_{\lambda}$ is the number of
reflections of $W_0^{\lambda}.$ For  $\lambda,\lambda'\in \Lambda^+$
we define
$$a_{\lambda,\lambda'}=\frac{v^{\nu_{\lambda'}}}{\pi_{\lambda'}}\sum_{w\in W_0}
(-1)^{l(w)}\Phi(\lambda'+\rho-w(\lambda+\rho)).$$

 \noindent {\bf Remark 6.1} (1) For any element in $\Lambda$, we will
use the same notation when it is regarded as an element in
$\widetilde{W}$. For two elements in $\L$, the operation between
them will be written additively if they are regarded as elements in
$\Lambda$ and will be written  multiplicatively when they are
regarded as elements in $\widetilde{W}$.\\
(2) In [L3], Lusztig showed that there is a 1-1 correspondence
between $\Lambda^+_r=\Lambda^+\bigcap\Lambda_r$ and the set of
$W_0-W_0$ double cosets in the affine Weyl group $W$. When we
consider the extended affine Weyl group $\widetilde{W}$, there is a
1-1 correspondence between $\Lambda^+$ and the set of $W_0-W_0$
double cosets in $\widetilde{W}$. So we have the similar results in
$\widetilde{W}$ for $a_{\l,\l'}$ and $b_{\l,\l'}$.

There is a 1-1 correspondence between $\Lambda^+$ and the set of
$W_0-W_0$ double cosets in $\widetilde{W}$ (an element $\lambda$ of
$\Lambda^+$ corresponds to the unique double coset $W_0\lambda W_0$
containing it: see [L5, \S2]). For each $\lambda\in \Lambda^+,$
there is a unique element $m_\lambda$ of minimal length and a unique
element $M_\lambda$ of maximal length in $W_0\lambda W_0.$ We have
$\lambda\leq\lambda'$ (i.e. $\l'-\l\in \mathbb NR^+$) if and only if
$M_\lambda\leq M_{\lambda'}$ (Bruhat order) for $\lambda,\lambda'\in
\Lambda^+$.

Set $v=q^{\frac{1}{2}},\ p_{u,w}=v^{l(u)-l(w)}P_{u,w}(v^{2})\in
\mathbb{Z} [v^{-1}]$ for $u\leq w\in W$ and $p_{u,w}=0$ for all
other $u,w$ in $W.$ For two elements $\l,\ \l'$ in $\L^+$, set
$$b_{\lambda,\lambda'}=\sum_{z\in W_0\lambda
W_0}(-v)^{l(m_\lambda)-l(z)}p_{z,m_{\lambda'}}.$$

\noindent {\bf Lemma 6.2 }([L3, \S5]) We have \\
 \noindent (1)
$a_{\lambda,\lambda'}$ is zero unless $\lambda\leq\lambda'$ and is
equal to 1 when $\lambda=\lambda';$ when $\lambda<\lambda'$ it
belongs to $v^{-1}\mathbb{Z}[v^{-1}].$\\
\noindent (2) $b_{\lambda,\lambda'}$ is zero unless
$\lambda\leq\lambda'$ and is equal to 1 when $\lambda=\lambda';$
when $\lambda<\lambda'$ it belongs to $v^{-1}\mathbb{Z} [v^{-1}].$
Also we have
$\mu(m_\lambda,m_{\lambda'})=\textrm{Res}_{v=0}(b_{\lambda,\lambda'}),$
where $\textrm{Res}_{v=0}(f)\in \mathbb{Z} $ denotes the coefficient
of $v^{-1}$ in $f\in {\mathcal {A}}=\mathbb
Z[v,v^{-1}].$\hfill$\Box$\

Let $\bar{}:\ {\mathcal {A}}\rightarrow{\mathcal {A}} $ be the ring
involution such that $\bar{v}=v^{-1}.$ The following lemma of
Lusztig gives a way to compute $b_{\l,\l'}$ inductively.

\noindent {\bf Lemma 6.3} ([L3, Proposition 7]) For any
$\lambda,\lambda''\in \Lambda^+,$ we have
$$\sum_{\lambda'\in
\Lambda_r^+}a_{\lambda,\lambda'}(-1)^{l(m_{\lambda'})-l(M_{\lambda'})}
\pi_{\lambda'}\bar{b}_{\lambda',\lambda''}=\sum_{\lambda'\in
\Lambda_r^+}\bar{a}_{\lambda,\lambda'}(-1)^{l(m_{\lambda'})-l(M_{\lambda'})}
\pi_{\lambda'}b_{\lambda',\lambda''}.$$\hfill$\Box$\

\def\L{\Lambda}
\def\l{\lambda}
\def\a{\alpha}
\def\b{\beta}
\def\r{\gamma}
\def\n{\eta}

 Then the set of positive roots in $R$ is
$$R^+=\{\a,\b,\a+\b\}.$$
We have $$\rho=\frac 12\alpha_{R^+}=\a+\b=x+y,$$ and
$$s(x)=-x+y,s(y)=y,t(x)=x,t(y)=x-y.$$

Define
\begin{eqnarray*}
X_1&=&\{mx\mid m\geq1\}\\
X_2&=&\{ny\mid n\geq1\}\\
Y_1&=&\{mx+y\mid m\geq1\}\\
Y_2&=&\{x+ny\mid n\geq1\}\\
Z_1&=&\{m(x+y)+ny\mid m\geq2,n\geq0\}\\
Z_2&=&\{n(x+y)+mx\mid n\geq2,m>0\}.
\end{eqnarray*}
Set $X=X_1\bigsqcup X_2,\ Y=Y_1\bigsqcup Y_2$, and $Z=Z_1\bigsqcup
Z_2$. We have $\Lambda_r^{+}=X\bigsqcup Y\bigsqcup Z$.

The following result is needed in calculating $a_{\l,\l'}$ and
$b_{\l',\l}$.

\noindent {\bf Proposition 6.4} Let $\lambda\in\Lambda^{+}$. We have\\
(1) $W_{0}^{\lambda}=\left\{\begin{array}{lll}
\{e,t\},& \ \ \textrm{if} \ \lambda\in X_1,\\
\{e,s\},& \ \ \textrm{if} \ \lambda\in X_2,\\
0,& \ \ \textrm{if} \ \lambda\in Y\bigsqcup Z.\end{array} \right.$\\
 (2)
$\nu_{\lambda}=\left\{\begin{array}{ll}
1,& \ \ \textrm{if} \ \lambda\in X,\\
0,& \ \ \textrm{if} \ \lambda\in Y\bigsqcup Z.\end{array} \right.$\\
(3) $\pi_{\lambda}=\left\{\begin{array}{ll}
{v+v^{-1}},& \ \ \textrm{if} \ \lambda\in X,\\
1,& \ \ \textrm{if} \ \lambda\in Y\bigsqcup Z.\end{array} \right.$\\
(4) For $\l\in X_1$ and $k\geq0$, we
have$${m_{\l}}=\left\{\begin{array}{ll}
(rst)^{2k+1}r,& \ \ \textrm{if} \ \lambda=(3k+3)x,\\
(rst)^{2k}\omega,& \ \ \textrm{if} \ \lambda=(3k+1)x,\\
(rst)^{2k}rs\omega^2,& \ \ \textrm{if} \ \lambda=(3k+2)x.
\end{array} \right.$$\\
(5) For $\l\in X_2$ and $k\geq0$, we
have$${m_{\l}}=\left\{\begin{array}{ll}
(rts)^{2k+1}r,& \ \ \textrm{if} \ \lambda=(3k+3)y,\\
(rts)^{2k}\omega^2,& \ \ \textrm{if} \ \lambda=(3k+1)y,\\
(rts)^{2k}rt\omega,& \ \ \textrm{if} \ \lambda=(3k+2)y.
\end{array} \right.$$\\
(6) For $\l\in Y_1$ and $k\geq0$, we
have$${m_{\l}}=\left\{\begin{array}{ll}
(rst)^{2k}r,& \ \ \textrm{if} \ \lambda=(3k+1)x+y,\\
(rst)^{2k+1}\omega,& \ \ \textrm{if} \ \lambda=(3k+2)x+y,\\
(rst)^{2k+1}rs\omega^2,& \ \ \textrm{if} \ \lambda=(3k+3)x+y.
\end{array} \right.$$\\
(7) For $\l\in Y_2$ and $k\geq0$, we
have$${m_{\l}}=\left\{\begin{array}{ll}
(rts)^{2k}r,& \ \ \textrm{if} \ \lambda=x+(3k+1)y,\\
(rts)^{2k+1}\omega^2,& \ \ \textrm{if} \ \lambda=x+(3k+2)y,\\
(rts)^{2k+1}rt\omega,& \ \ \textrm{if} \ \lambda=x+(3k+3)y.
\end{array} \right.$$\\
(8) For $\l\in Z_1$ and $m\geq2,k\geq0$, we
have$${m_{\l}}=\left\{\begin{array}{ll}
r(stsr)^{m-1}(tsr)^{2k},& \ \ \textrm{if} \ \lambda=mx+(m+3k)y,\\
r(stsr)^{m-1}(tsr)^{2k}ts\omega^2,& \ \ \textrm{if} \ \lambda=mx+(m+3k+1)y,\\
r(stsr)^{m-1}(tsr)^{2k}t\omega,& \ \ \textrm{if} \
\lambda=mx+(m+3k+2)y.
\end{array} \right.$$\\
(9) For $\l\in Z_2$ and $n\geq2,k>0$, we
have$${m_{\l}}=\left\{\begin{array}{ll}
r(stsr)^{n-1}(str)^{2k},& \ \ \textrm{if} \ \lambda=(3k+n)x+ny,\\
r(stsr)^{n-1}(str)^{2k}st\omega,& \ \ \textrm{if} \ \lambda=(3k+1+n)x+ny,\\
r(stsr)^{n-1}(str)^{2k}s\omega^2,& \ \ \textrm{if} \
\lambda=(3k+2+n)x+ny.
\end{array} \right.$$\\
(10) $(-1)^{l(m_{\lambda})-l(M_{\lambda})}=\left\{\begin{array}{ll}
-1,& \ \ \textrm{if} \ \lambda\in X,\\
1,& \ \ \textrm{if} \ \lambda\in Y\bigsqcup Z.\end{array} \right.$

\noindent {\bf Proof.} We have that $s(x)=x-\a=y-x,s(y)=y,t(x)=x$,
and $t(y)=y-\b=x-y$. Assume that $\l=mx+ny$, $m,n\in\mathbb{N}$.
Then we get that $s(\l)=-mx+(m+n)y$, $t(\l)=(m+n)x-ny$,
$st(\l)=-(m+n)x+my$, $ts(\l)=nx-(m+n)y$, and $sts(\l)=-\l$. Thus one
can check (1)-(3) easily.

We now prove (4). First, the corresponding elements of $x,2x,3x$ in
$\widetilde{W}$ are $\omega tr, (\omega tr)^2=\omega^2rstr$ and
$(\omega tr)^3=strstr$, respectively. Hence, the corresponding
elements of $(3k+1)x,(3k+2)x,3kx$ in $\widetilde{W}$ are $(\omega
tr)^{3k+1}=(strstr)^k\omega tr, (\omega
tr)^{3k+2}=(strstr)^k\omega^2rstr$ and $(\omega
tr)^{3k}=(strstr)^k$, respectively.

Thus, the minimal element $m_{\l}$ in  the double costs $W_0\l W_0$
is$${m_{\l}}=\left\{\begin{array}{ll}
(rst)^{2k+1}r,& \ \ \textrm{if} \ \lambda=(3k+3)x,\\
(rst)^{2k}\omega,& \ \ \textrm{if} \ \lambda=(3k+1)x,\\
(rst)^{2k}rs\omega^2,& \ \ \textrm{if} \ \lambda=(3k+2)x.
\end{array} \right.$$
Moreover,  the maximal element $M_{\l}$ in  the double costs $W_0\l
W_0$ is$${M_{\l}}=\left\{\begin{array}{ll}
sts(rst)^{2k+1}rst,& \ \ \textrm{if} \ \lambda=(3k+3)x,\\
sts(rst)^{2k}\omega st,& \ \ \textrm{if} \ \lambda=(3k+1)x,\\
sts(rst)^{2k}rs\omega^2 st,& \ \ \textrm{if} \ \lambda=(3k+2)x.
\end{array} \right.$$
We get that $(-1)^{l(m_{\lambda})-l(M_{\lambda})}=-1$.

With the similar method,we can prove (5-9) and can get
(10).\hfill$\Box$\

Simple computations lead to the following identities:\\
$\Phi(0)=1,$\\
$\Phi(\alpha)=\Phi(\beta)=-v^{-2},$\\
$\Phi(\alpha+\beta)=v^{-4}-v^{-2},$\\
$\Phi(2\alpha+\beta)=\Phi(\alpha+2\beta)=v^{-4},$\\
$\Phi(2\alpha+2\beta)=-v^{-6},$\\
and $\Phi(\lambda)=0$ for all other $\lambda$ in $\Lambda_r$.

Using the above formulas we can compute $a_{\l,\l'}$, which are
needed in determining $b_{\l,\l''}$ in the next section.

\noindent {\bf Proposition 6.5} Let $0<\l<\l'\in\Lambda^{+}$. Assume
that
$\lambda=mx+ny$, $m,n\in\mathbb{N}$. \\
(1) If $\l'-\l$ is not in $\{\a,\b,\a+\b,2\a+\b,\a+2\b,2\a+2\b\}$,
then
 $a_{\l,\l'}=0$.\\
(2) If $\l'-\l=\a$ or $\b$, then $a_{\l,\l'}=-v^{-2}$.\\
(3) If $\l'-\l=\a+\b$, then
$a_{\lambda,\lambda'}=\left\{\begin{array}{ll}
-v^{-2},& \ \ \textrm{when} \ mn=0,\\
v^{-4}-v^{-2},& \ \ \textrm{otherwise}.\end{array} \right.$\\
(4) If $\l'-\l=\a+2\b$ or $2\a+\b$, then $a_{\l,\l'}=v^{-4}$.\\
(5) If $\l'-\l=2\a+2\b$, then $a_{\l,\l'}=-v^{-6}$.

\noindent {\bf Proof.} Since $\rho=\a+\b$, we get that
$\rho-s(\rho)=\a,$ $\rho-t(\rho)=\b,$ $\rho-st(\rho)=2\a+\b,$
$\rho-ts(\rho)=\a+2\b,$ and $\rho-sts(\rho)=2\a+2\b$.

We also have that $e(\l)=\l=mx+ny$, $s(\l)=-mx+(m+n)y,$
$t(\l)=(m+n)x-ny,$ $st(\l)=-(m+n)x+my,$ $ts(\l)=nx-(m+n)y,$ and
$sts(\l)=-mx-ny$.

(a) Assume that $\l'=m'x+n'y,m',n'\in\mathbb{N}$. Since
$x=\frac{2}{3}\alpha+\frac{1}{3}\beta,
y=\frac{1}{3}\alpha+\frac{2}{3}\beta$, we  have
$\l'+\rho-st(\l+\rho)=\l'-st(\l)+2\a+\b=\frac{1}{3}[(2m'+n'+m+2n)
\a+(m'+2n'+n-m)\b]+2\a+\b$. Since $2m'+n'+m+2n>0$, the coefficient
of $\a$ in $\l'+\rho-st(\l+\rho)$ is $>2$, which implies that
$\Phi(\l'+\rho-st(\l+\rho))=0$. Similarly, we can show that
$\Phi(\l'+\rho-ts(\l+\rho))=\Phi(\l'+\rho-sts(\l+\rho))=0$.

Suppose that $\l'-\l=i\a+j\b,i,j\in\mathbb{N}$. We have that
$\l'+\rho-s(\l+\rho)=(m+i+1)\a+j\b$ and
$\l'+\rho-t(\l+\rho)=i\a+(n+j+1)\b$.

(b) If  $\l'-\l=\a$, we get $\l'=(m+2)x+(n-1)y, n\geq1$. We first
have that $\Phi(\l'-\l)=-v^{-2}$. Also we have
$\l'+\rho-s(\l+\rho)=(m+2)\a$, which implies that
$\Phi(\l'+\rho-s(\l+\rho))=0$. Since
$\l'+\rho-t(\l+\rho)=\a+(n+1)\b$, we get
$$\Phi(\l'+\rho-t(\l+\rho))=\left\{\begin{array}{ll}
v^{-4},& \ \ \textrm{if}\ n=1,\\
0,& \ \ \textrm{if} \ n\geq2.\end{array} \right.$$When $n=1$, we
have $\pi_{\l'}=v+v^{-1}$ and $\nu_{\l'}=1$ by Proposition 2.4
(2-3). Thus we get
$a_{\l,\l'}=\frac{v}{v+v^{-1}}(-v^{-2}-v^{-4})=-v^{-2}$.

 When $n\geq2$, we have $\pi_{\l'}=1$ and $\nu_{\l'}=0$. Thus we get
$a_{\l,\l'}=-v^{-2}$.

(c) If  $\l'-\l=\b$, we get $\l'=(m-1)x+(n+2)y, m\geq1$. We first
have that $\Phi(\l'-\l)=-v^{-2}$. Also we have
$\l'+\rho-t(\l+\rho)=(n+2)\b$, which implies that
$\Phi(\l'+\rho-t(\l+\rho))=0$. Since
$\l'+\rho-s(\l+\rho)=(m+1)\a+\b$, we get
$$\Phi(\l'+\rho-s(\l+\rho))=\left\{\begin{array}{ll}
v^{-4},& \ \ \textrm{if}\ m=1,\\
0,& \ \ \textrm{if} \ m\geq2.\end{array} \right.$$When $m=1$, we
have $\pi_{\l'}=v+v^{-1}$ and $\nu_{\l'}=1$. Thus we get
$a_{\l,\l'}=\frac{v}{v+v^{-1}}(-v^{-2}-v^{-4})=-v^{-2}$.

 When $m\geq2$, we have $\pi_{\l'}=1$ and $\nu_{\l'}=0$. Thus we get
$a_{\l,\l'}=-v^{-2}$.

By (b) and (c), we prove (2)

(d) If $\l'-\l=\a+\b$, we get $\l'=(m+1)x+(n+1)y$. We first have
that $\Phi(\l'-\l)=v^{-4}-v^{-2}$. Since
$\l'+\rho-s(\l+\rho)=(m+2)\a+\b$, we get
$$\Phi(\l'+\rho-s(\l+\rho))=\left\{\begin{array}{ll}
v^{-4},& \ \ \textrm{if}\ m=0,\\
0,& \ \ \textrm{if} \ m\geq1.\end{array} \right.$$ Since
$\l'+\rho-t(\l+\rho)=\a+(n+2)\b$, we get
$$\Phi(\l'+\rho-s(\l+\rho))=\left\{\begin{array}{ll}
v^{-4},& \ \ \textrm{if}\ n=0,\\
0,& \ \ \textrm{if} \ n\geq1.\end{array} \right.$$ Since
$\pi_{\l'}=1,\nu_{\l'}=0$ and $m+n\geq1$, we get
$$a_{\lambda,\lambda'}=\left\{\begin{array}{ll}
-v^{-2},& \ \ \textrm{when} \ mn=0,\\
v^{-4}-v^{-2},& \ \ \textrm{otherwise}.\end{array} \right.$$Thus (3)
holds.

(e) If $\l'-\l=2\a+\b$, we get $\l'=(m+3)x+ny$. We first have that
$\Phi(\l'-\l)=v^{-4}$. Also we have $\l'+\rho-s(\l+\rho)=(m+3)\a+\b$
and $m+3\geq3$, which implies that $\Phi(\l'+\rho-s(\l+\rho))=0$.
Since $\l'+\rho-t(\l+\rho)=2\a+(n+2)\b$, we get
$$\Phi(\l'+\rho-t(\l+\rho))=\left\{\begin{array}{ll}
-v^{-6},& \ \ \textrm{if}\ n=0,\\
0,& \ \ \textrm{if} \ n\geq1.\end{array} \right.$$When $n=1$, we see
that $\pi_{\l'}=v+v^{-1}$ and $\nu_{\l'}=1$. Thus we get
$a_{\l,\l'}=\frac{v}{v+v^{-1}}(v^{-4}+v^{-6})=v^{-4}$.

When $n\geq1$, we have $\pi_{\l'}=1$ and $\nu_{\l'}=0$, which
implies that $a_{\l,\l'}=v^{-4}$.

(f) If $\l'-\l=\a+2\b$, we get $\l'=mx+(n+3)y$. We first have that
$\Phi(\l'-\l)=v^{-4}$. Also we have $\l'+\rho-t(\l+\rho)=\a+(n+3)\b$
and $n+3\geq3$, which implies that $\Phi(\l'+\rho-t(\l+\rho))=0$.
Since $\l'+\rho-s(\l+\rho)=(m+2)\a+2\b$, we get
$$\Phi(\l'+\rho-s(\l+\rho))=\left\{\begin{array}{ll}
-v^{-6},& \ \ \textrm{if}\ m=0,\\
0,& \ \ \textrm{if} \ m\geq1.\end{array} \right.$$When $m=1$, we see
that $\pi_{\l'}=v+v^{-1}$ and $\nu_{\l'}=1$. Thus we get
$a_{\l,\l'}=\frac{v}{v+v^{-1}}(v^{-4}+v^{-6})=v^{-4}$.

When $m\geq1$, we have $\pi_{\l'}=1$ and $\nu_{\l'}=0$, which
implies that $a_{\l,\l'}=v^{-4}$.

By (e) and (f), we prove (4).

(g) If $\l'-\l=2\a+2\b$, we get $\l'=(m+2)x+(n+2)y$. We first have
that $\Phi(\l'-\l)=-v^{-6}$. Since $\l'+\rho-s(\l+\rho)=(m+3)\a+2\b$
and $\l'+\rho-s(\l+\rho)=2\a+(n+3)\b$, we get
$\Phi(\l'+\rho-s(\l+\rho))=\Phi(\l'+\rho-t(\l+\rho))=0$.

By $\pi_{\l'}=1$ and $\nu_{\l'}=0$, we get $a_{\l,\l'}=-v^{-6}$.
Thus (5) holds.

(h) If $\l'-\l=i\a+j\b$, $i,j\in\mathbb{N}$ for three cases: (a)
$i=2$ and $j=0$; (b) $i=0$ and $j=2$; (c) $i$ or $j\geq3$. For these
three cases, we always have
$\Phi(\l'-\l)=\Phi(\l'+\rho-s(\l+\rho))=\Phi(\l'+\rho-t(\l+\rho))=0$,
which implies that $a_{\l,\l'}=0$. Thus (1) holds. \hfill$\Box$\

We have the set of dominant weights is
\begin{eqnarray*}
\L^+&=&\{mx+ny,\
m,n\in\mathbb N\}\\
&=&\{\frac{1}{3}(i\a+j\b)\mid i,j\in\mathbb{N},\frac{i}{2}\leq
j\leq2i,3\mid2i-j,3\mid2j-i \}.
\end{eqnarray*}

Now assume that $\l''\in\Lambda^+$. We will compute $b_{\l,\l''}$
for any $0<\l\leq\l''$ in $\Lambda^+$ by Lemma 6.3. For any
$\l<\l''\in \Lambda^+$, we denote $$c_{\l,\l''}=\sum_{\lambda'\in
\Lambda_r^+}a_{\lambda,\lambda'}(-1)^{l(m_{\lambda'})-l(M_{\lambda'})}
\pi_{\lambda'}\bar{b}_{\lambda',\lambda''},$$ the left hand side in
the equation in Lemma 2.3. We set $\xi=v+v^{-1}$.

\noindent {\bf Proposition 6.6} Let
$\l''=ny=\frac{n}{3}\a+\frac{2n}{3}\b\in X$ with $n\geq4$. Then for
any $0<\l<\l''$ in $\Lambda^+$, we have
$${b_{\l,\l''}}=\left\{\begin{array}{ll}
{v^{-1}+v^{-3}}, & \ \ \textrm{if} \;  \l=\l''-\b ,\\
{v^{-4}}, & \ \ \textrm{if} \;  \l=\l''-(\a+2\b) ,\\
{0} ,   & \ \ \textrm{otherwise}.\end{array} \right.$$

\noindent {\bf Proof.} We will compute all $b_{\l,\l''}$ by using
Lemma 6.3. Lemma 6.4 (2-3) and Lemma 6.5 will be used in the
following computations.

Set $\l_j=\frac{n}{3}\a+\frac{2n}{3}\b-j\b$,
$\r_j=(\frac{n}{3}-1)\a+(\frac{2n}{3}-2)\b-j\b$,
$\n_j=(\frac{n}{3}-2)\a+(\frac{2n}{3}-4)\b-j\b\in\Lambda^+$, with
$j\geq0$. We see that $\l''=\l_0$.

\noindent(1) We first compute $b_{\l_1,\l_0}$. Since $\l_1$ is not
in $X$ and $\l_0-\l_1=\b$, we have
\begin{eqnarray*}
c_{\l_1,\l_0}&=&a_{\l_1,\l_1}(-1)^{l(m_{\l_1})-l(M_{\l_1})}\pi_{\l_1}
\overline{b}_{\l_1,\l_0}+a_{\l_1,\l_0}(-1)^{l(m_{\l_0})-l(M_{\l_0})}\pi_{\l_0}
\overline{b}_{\l_0,\l_0}\\
&=&\overline{b}_{\l_1,\l_0}-v^{-2}(-\xi)\\
&=&\overline{b}_{\l_1,\l_0}+v^{-1}+v^{-3}.
\end{eqnarray*}
By Lemma 6.3, we have
$\overline{b}_{\l_1,\l_0}+v^{-1}+v^{-3}=b_{\l_1,\l_0}+v+v^{3}$.
Since $b_{\l_1,\l_0}\in v^{-1}\mathbb{Z} [v^{-1}]$, we get
$b_{\l_1,\l_0}=v^{-1}+v^{-3}$.

For $\l_2$, we have $\l_1-\l_2=\b$ and $\l_0-\l_2=2\b$. Thus we get
\begin{eqnarray*}
c_{\l_2,\l_0}&=&(-1)^{l(m_{\l_2})-l(M_{\l_2})}\pi_{\l_2}
\overline{b}_{\l_2,\l_0}+a_{\l_2,\l_1}(-1)^{l(m_{\l_1})-l(M_{\l_1})}\pi_{\l_1}
\overline{b}_{\l_1,\l_0}\\
&=&(-1)^{l(m_{\l_2})-l(M_{\l_2})}\pi_{\l_2}\overline{b}_{\l_2,\l_0}+(-v^{-2})(v+v^{3})\\
&=&(-1)^{l(m_{\l_2})-l(M_{\l_2})}\pi_{\l_2}\overline{b}_{\l_2,\l_0}-\xi.
\end{eqnarray*}
Since $\overline{\xi}=\xi$, we get
$\overline{b}_{\l_2,\l_0}=b_{\l_2,\l_0}$ by Lemma 6.3, which implies
that $b_{\l_2,\l_0}=0$.

Assume that $b_{\l_k,\l_0}=0$ for any $2\leq k\leq j-1$. Now we
prove that $b_{\l_j,\l_0}=0$ for $j\geq3$. Since
$a_{\l_j,\l_j}=1,a_{\l_j,\l_1}=a_{\l_j,\l_0}=0$, we get
$c_{\l_j,\l_0}=(-1)^{l(m_{\l_j})-l(M_{\l_j})}\pi_{\l_j}
\overline{b}_{\l_j,\l_0}$. By Lemma 6.3, we get
$\overline{b}_{\l_j,\l_0}=b_{\l_j,\l_0}$. Thus $b_{\l_j,\l_0}=0$ for
$j\geq3$.

\noindent(2) Now we compute $b_{\r_0,\l_0}$. Since $\r_0\in X$,
$\l_1-\r_0=\a+\b$ and $\l_0-\r_0=\a+2\b$, we have
\begin{eqnarray*}
c_{\r_0,\l_0}&=&(-1)^{l(m_{\r_0})-l(M_{\r_0})}\pi_{\r_0}
\overline{b}_{\r_0,\l_0}+a_{\r_0,\l_1}(-1)^{l(m_{\l_1})-l(M_{\l_1})}\pi_{\l_1}
\overline{b}_{\l_1,\l_0}\\
&&+a_{\r_0,\l_0}(-1)^{l(m_{\l_0})-l(M_{\l_0})}\pi_{\l_0}
\overline{b}_{\l_0,\l_0}\\
&=&-\xi\overline{b}_{\r_0,\l_0}+(-v^{-2})(v+v^{3})+v^{-4}(-\xi)\\
&=&-\xi\overline{b}_{\r_0,\l_0}-\xi+v^{-4}(-\xi).
\end{eqnarray*}
By Lemma 6.3, we can get $b_{\r_0,\l_0}=v^{-4}$.

For $\r_1$, we have $\r_0-\r_1=\b$, $\l_1-\r_1=\a+2\b$, and
$\l_0-\r_1=\a+3\b$. Thus we get
\begin{eqnarray*}
c_{\r_1,\l_0}&=&(-1)^{l(m_{\r_1})-l(M_{\r_1})}\pi_{\r_1}
\overline{b}_{\r_1,\l_0}+a_{\r_1,\r_0}(-1)^{l(m_{\r_0})-l(M_{\r_0})}\pi_{\r_0}
\overline{b}_{\r_0,\l_0}\\
&&+a_{\r_1,\l_1}(-1)^{l(m_{\l_1})-l(M_{\l_1})}\pi_{\l_1}
\overline{b}_{\l_1,\l_0}\\
&=&(-1)^{l(m_{\r_1})-l(M_{\r_1})}\pi_{\r_1}\overline{b}_{\r_1,\l_0}
+(-v^{-2})(-\xi)v^{4}+v^{-4}(v+v^{3})\\
&=&(-1)^{l(m_{\r_1})-l(M_{\r_1})}\pi_{\r_1}\overline{b}_{\r_1,\l_0}
+\xi+v^{3}+v^{-3}.
\end{eqnarray*}
By Lemma 6.3, we can get $b_{\r_1,\l_0}=0$.

For $\r_j$ with $j\geq2$, we have $a_{\r_j,\l}=0$ for
$\l=\r_0,\l_1,\l_0$. Thus we can get $b_{\r_j,\l_0}=0$ by induction
on the order $\leq$ in $\Lambda^+$ by Lemma 6.3.

\noindent(3) For $\n_0$, we have $\r_0-\n_0=\a+2\b$ and
$a_{\n_0,\l}=0$ for $\l=\l_1,\l_0$. Thus we can get
\begin{eqnarray*}
c_{\n_0,\l_0}&=&(-1)^{l(m_{\n_0})-l(M_{\n_0})}\pi_{\n_0}
\overline{b}_{\n_0,\l_0}+a_{\n_0,\r_0}(-1)^{l(m_{\r_0})-l(M_{\r_0})}\pi_{\r_0}
\overline{b}_{\r_0,\l_0}\\
&=&-\xi\overline{b}_{\n_0,\l_0}+v^{-4}(-\xi)v^{4}\\
&=&-\xi\overline{b}_{\n_0,\l_0}-\xi.
\end{eqnarray*}
By Lemma 6.3, we can get $b_{\n_0,\l_0}=0$.

For $\n_j$ with $j\geq1$, we have $a_{\n_j,\l}=0$ for
$\l=\r_0,\l_1,\l_0$. Thus we can get $b_{\n_j,\l_0}=0$ by induction
on the order $\leq$ in $\Lambda^+$ by Lemma 6.3.

\noindent(4) If $0<\zeta=\l_0-k(\a+2\b)-j\b\in\Lambda^+$ with
$k\geq3$, $j\geq0$. Since $a_{\zeta,\l}=0$ for $\l=\r_0,\l_1,\l_0$,
 we can get $b_{\zeta,\l_0}=0$ by Lemma 6.3. We finish the
proof. \hfill$\Box$\

\noindent {\bf Proposition 6.7} Suppose that
$\l''=x+my=\frac{m+2}{3}\a+(\frac{2m+4}{3}-1)\b\in\Lambda^+$ with
$m\geq4$. Then for any $0<\l<\l''$ in $\Lambda^+$, we have
$${b_{\l,\l''}}=\left\{\begin{array}{ll}
{-v^{-2}}, & \ \ \textrm{if} \;  \l=\l''-\b,\\
{v^{-1}}, & \ \ \textrm{if} \;  \l=\l''-(\a+\b),\\
{v^{-4}}, & \ \ \textrm{if} \;  \l=\l''-(\a+2\b),\\
{0} ,   & \ \ \textrm{otherwise} .\end{array} \right.$$

\noindent {\bf Proof.} We set
$\l_j=\frac{m+2}{3}\a+\frac{2m+4}{3}\b-j\b\in\Lambda^+$,
$\r_j=(\frac{m+2}{3}-1)\a+(\frac{2m+4}{3}-2)\b-j\b\in\Lambda^+$, and
$\n_j=(\frac{m+2}{3}-2)\a+(\frac{2m+4}{3}-4)\b-j\b\in\Lambda^+$,
with $j\geq0$. We see that $\l''=\l_1$.

\noindent(1) First we compute $b_{\l_2,\l''}$. Since $\l''-\l_2=\b$
and $\l_2,\l''$ are not in $X$, we get
\begin{eqnarray*}
c_{\l_2,\l''}&=&(-1)^{l(m_{\l_2})-l(M_{\l_2})}\pi_{\l_2}
\overline{b}_{\l_2,\l''}+a_{\l_2,\l''}(-1)^{l(m_{\l''})-l(M_{\l''})}\pi_{\l''}
\overline{b}_{\l'',\l''}\\
&=&\overline{b}_{\l_2,\l''}-v^{-2}.
\end{eqnarray*}
By Lemma 6.3, we get $b_{\l_2,\l''}=-v^{-2}$.

For $\l_3$, we have $\l_2-\l_3=\b$ and $\l''-\l_3=2\b$. Thus we get
\begin{eqnarray*}
c_{\l_3,\l''}&=&(-1)^{l(m_{\l_3})-l(M_{\l_3})}\pi_{\l_3}
\overline{b}_{\l_3,\l''}+a_{\l_3,\l_2}(-1)^{l(m_{\l_2})-l(M_{\l_2})}\pi_{\l_2}
\overline{b}_{\l_2,\l''}\\
&=&(-1)^{l(m_{\l_3})-l(M_{\l_3})}\pi_{\l_3}
\overline{b}_{\l_3,\l''}-v^{-2}(-v^2).
\end{eqnarray*}
By Lemma 6.3, we get $b_{\l_3,\l''}=0$.

For $\l_j$ with $j\geq4$,  we have $a_{\l_j,\l_2}=a_{\l_j,\l''}=0$.
Thus we can show that $b_{\l_j,\l''}=0$ by induction on the order
$\leq$.

\noindent(2) We now compute $b_{\r_0,\l''}$. Since $\r_0\in X$,
$\l_2-\r_0=\a$, and $\l''-\r_0=\a+\b$, we get
\begin{eqnarray*}
c_{\r_0,\l''}&=&(-1)^{l(m_{\r_0})-l(M_{\r_0})}\pi_{\r_0}
\overline{b}_{\r_0,\l''}+a_{\r_0,\l_2}(-1)^{l(m_{\l_2})-
l(M_{\l_2})}\pi_{\l_2}\overline{b}_{\l_2,\l''}\\
&&+a_{\r_0,\l''}(-1)^{l(m_{\l''})-l(M_{\l''})}\pi_{\l''}\\
&=&-\xi\overline{b}_{\r_0,\l''}-v^{-2}(-v^2)-v^{-2}\\
&=&-\xi\overline{b}_{\r_0,\l''}+1-v^{-2}.
\end{eqnarray*}
 By Lemma 6.3, we get $b_{\r_0,\l''}=v^{-1}$.

For $\r_1$, we have $\r_1$ is not in $X$, $\r_0-\r_1=\b$,
$\l_2-\r_1=\a+\b$, and $\l''-\r_1=\a+2\b$. Thus we get
\begin{eqnarray*}
c_{\r_1,\l''}&=&(-1)^{l(m_{\r_1})-l(M_{\r_1})}\pi_{\r_1}
\overline{b}_{\r_1,\l''}+a_{\r_1,\r_0}(-1)^{l(m_{\r_0})-l(M_{\r_0})}\pi_{\r_0}
\overline{b}_{\r_0,\l''}\\
&&+a_{\r_1,\l_2}(-1)^{l(m_{\l_2})-
l(M_{\l_2})}\pi_{\l_2}\overline{b}_{\l_2,\l''}+a_{\r_1,\l''}
(-1)^{l(m_{\l''})-l(M_{\l''})}\pi_{\l''}\\
&=&\overline{b}_{\r_1,\l''}-v^{-2}(-\xi)v+(v^{-4}-v^{-2})(-v^{2})+v^{-4}\\
&=&\overline{b}_{\r_1,\l''}+v^{-4}.
\end{eqnarray*}
By Lemma 6.3, we get $b_{\r_1,\l''}=v^{-4}$.

For $\r_2$, we have  $\r_1-\r_2=\b$, $\r_0-\r_2=2\b$,
$\l_2-\r_2=\a+2\b$, and $\l''-\r_2=\a+3\b$. Thus we get
\begin{eqnarray*}
c_{\r_2,\l''}&=&(-1)^{l(m_{\r_2})-l(M_{\r_2})}\pi_{\r_2}
\overline{b}_{\r_2,\l''}+a_{\r_2,\r_1}(-1)^{l(m_{\r_1})-l(M_{\r_1})}\pi_{\r_1}
\overline{b}_{\r_1,\l''}\\
&&+a_{\r_2,\l_2}(-1)^{l(m_{\l_2})-
l(M_{\l_2})}\pi_{\l_2}\overline{b}_{\l_2,\l''}\\
&=&(-1)^{l(m_{\r_2})-l(M_{\r_2})}\pi_{\r_2}\overline{b}_{\r_1,\l_1}
-v^{-2}v^4+v^{-4}(-v^{2}).
\end{eqnarray*}
By Lemma 6.3, we get $b_{\r_2,\l''}=0$.

For $\r_j$ with $j\geq3$,  we have $a_{\r_j,\l}=0$ for
$\l=\r_1,\r_0,\l_2$ and $\l''$. Thus we can show that
$b_{\r_j,\l''}=0$ by induction on the order $\leq$.

\noindent(3) Now we compute $b_{\n_0,\l''}$. Since $\n_0\in X$,
$\r_1-\n_0=\a+\b$, $\r_0-\n_0=\a+2\b$, $\l_2-\n_0=2\a+2\b$, and
$\l''-\n_0=2\a+3\b$, we get
\begin{eqnarray*}
c_{\n_0,\l''}&=&(-\xi)\overline{b}_{\n_0,\l''}
+a_{\n_0,\r_1}(-1)^{l(m_{\r_1})-l(M_{\r_1})}\pi_{\r_1}
\overline{b}_{\r_1,\l''}\\
&&+a_{\n_0,\r_0}(-1)^{l(m_{\r_0})-l(M_{\r_0})}\pi_{\r_0}
\overline{b}_{\r_0,\l''}+a_{\n_0,\l_2}(-1)^{l(m_{\l_2})-
l(M_{\l_2})}\pi_{\l_2}\overline{b}_{\l_2,\l''}\\
&=&(-\xi)\overline{b}_{\n_0,\l''}-v^{-2}v^4+v^{-4}(-\xi)v-v^{-6}(-v^{2})\\
&=&(-\xi)\overline{b}_{\n_0,\l''}-v^{2}-v^{-2}-v^{-4}+v^{-4}.
\end{eqnarray*}
By Lemma 6.3, we get $b_{\n_0,\l''}=0$.

For $\n_1$, we have $\r_1-\n_1=\a+2\b$ and  $a_{\n_1,\l}=0$ for
$\l=\r_0,\l_2,\l''$. Thus we get
\begin{eqnarray*}
c_{\n_1,\l''}&=&(-1)^{l(m_{\n_1})-l(M_{\n_1})}\pi_{\n_1}\overline{b}_{\n_1,\l''}
+a_{\n_1,\r_1}(-1)^{l(m_{\r_1})-l(M_{\r_1})}\pi_{\r_1}
\overline{b}_{\r_1,\l''}\\
&=&(-1)^{l(m_{\n_1})-l(M_{\n_1})}\pi_{\n_1}\overline{b}_{\n_1,\l''}
+v^{-4}v^4.
\end{eqnarray*}
By Lemma 6.3, we get $b_{\n_1,\l''}=0$

For $\n_j$ with $j\geq2$, we have $a_{\n_j,\l}=0$ for
$\l=\r_1,\r_0,\l_2,\l''$. Thus we can show that $b_{\n_j,\l''}=0$ by
Lemma 6.3.

\noindent(4) If $\zeta=\l_0-k(\a+2\b)-j\b$ with $k\geq3,j\geq0$. We
have $a_{\zeta,\l}=0$ for $\l=\r_1,\r_0,\l_2,\l''$. Thus we can get
that $b_{\zeta,\l''}=0$ by Lemma 6.3.

We finish the proof. \hfill$\Box$\

\noindent {\bf Proposition 6.8} Let
$\l''=mx=\frac{2m}{3}\a+\frac{m}{3}\b\in X$ with $m\geq4$. Then for
$0<\l<\l''$ in $\Lambda^+$, we have
$${b_{\l,\l''}}=\left\{\begin{array}{ll}
{v^{-1}+v^{-3}}, & \ \ \textrm{if} \;  \l=\l''-\a ,\\
{v^{-4}}, & \ \ \textrm{if} \; \l=\l''-(2\a+\b),\\
{0} ,   & \ \ \textrm{otherwise} \;  .\end{array} \right.$$

\noindent {\bf Proof.} We set
$\r_j=(\frac{2m}{3}-1)\a+\frac{m}{3}\b+j\b$,
$\n_j=(\frac{2m}{3}-2)\a+(\frac{m}{3}-1)\b+j\b$,
$\theta_j=(\frac{2m}{3}-3)\a+(\frac{m}{3}-1)\b+j\b$, and
$\zeta_j=(\frac{2m}{3}-4)\a+(\frac{m}{3}-2)\b+j\b$ in $\Lambda^+$
with $j\geq0$.

(1) We see that $0<\r_j<\l''$ if and only if $j=0$.
 Now compute $b_{\r_0,\l''}$. Since $\r_0$ is
not in $X$, $\l''-\r_0=\a$, we have $a_{\r_0,\r_0}=1$ and
$a_{\r_0,\l''}=-v^{-2}$. Thus we get
\begin{eqnarray*}
c_{\r_0,\l''}&=&(-1)^{l(m_{\r_0})-l(M_{\r_0})}\pi_{\r_0}
\overline{b}_{\r_0,\l''}+a_{\r_0,\l''}(-1)^{l(m_{\l''})-l(M_{\l''})}\pi_{\l''}\\
&=&\overline{b}_{\r_0,\l''}-v^{-2}(-\xi)\\
&=&\overline{b}_{\r_0,\l''}+v^{-1}+v^{-3}.
\end{eqnarray*}
 By Lemma 6.3, we can get $b_{\r_0,\l''}=v^{-1}+v^{-3}$.

(2)  We see that $0<\n_j<\l''$ if and only if $j=0$ or $1$. First we
compute $b_{\n_1,\l''}$. Since $\r_0-\n_1=\a$ and $\l''-\n_1=2\a$,
we have $a_{\n_1,\r_0}=-v^{-2}$ and $a_{\n_1,\l''}=0$. Thus we get
\begin{eqnarray*}
c_{\n_1,\l''}&=&(-1)^{l(m_{\n_1})-l(M_{\n_1})}\pi_{\n_1}
\overline{b}_{\n_1,\l''}+a_{\n_1,\r_0}(-1)^{l(m_{\r_0})-l(M_{\r_0})}\pi_{\r_0}
\overline{b}_{\r_0,\l''}\\
&=&(-1)^{l(m_{\n_1})-l(M_{\n_1})}\pi_{\n_1}\overline{b}_{\n_1,\l''}-v^{-2}(v+v^{3})\\
&=&(-1)^{l(m_{\n_1})-l(M_{\n_1})}\pi_{\n_1}\overline{b}_{\n_1,\l''}-v^{-1}-v.
\end{eqnarray*}
By Lemma 6.3, we can get $b_{\n_1,\l''}=0$.

Now we compute $b_{\n_0,\l''}$. Since $\n_0\in X$,
$\r_0-\n_0=\a+\b$, and $\l''-\n_0=2\a+\b$, we have
$a_{\n_0,\r_0}=-v^{-2}$ and $a_{\n_0,\l''}=v^{-4}$. Thus we get
\begin{eqnarray*}
c_{\n_0,\l''}&=&(-1)^{l(m_{\n_0})-l(M_{\n_0})}\pi_{\n_0}
\overline{b}_{\n_0,\l''}+a_{\n_0,\r_0}(-1)^{l(m_{\r_0})-l(M_{\r_0})}\pi_{\r_0}
\overline{b}_{\r_0,\l''}\\
&&+a_{\n_0,\l''}(-1)^{l(m_{\l''})-l(M_{\l''})}\pi_{\l''}\\
&=&(-\xi)\overline{b}_{\n_0,\l''}-v^{-2}(v^{}+v^{3})+v^{-4}(-\xi).
\end{eqnarray*}
By Lemma 6.3, we can get $b_{\n_0,\l''}=v^{-4}$.

(3) We see that $0<\theta_j<\l''$ if and only if $j=0$ or $1$ and
$m\geq5$.

Obviously, we have $\theta_1\nless \n_0$, and
$a_{\theta_1,\r_0}=a_{\theta_1,\l''}=0$. Thus we have
$c_{\theta_1,\l''}=(-1)^{l(m_{\theta_1})-l(M_{\theta_1})}\pi_{\theta_1}
\overline{b}_{\theta_1,\l''}$, which implies that
$b_{\theta_1,\l''}=0$.

Since $\n_0-\theta_0=\a$,  $\r_0-\theta_0=2\a+\b$, and
$\l''-\theta_0=3\a+\b$, we have $a_{\theta_0,\n_0}=-v^{-2}$,
$a_{\theta_0,\r_0}=v^{-4}$, and $a_{\theta_0,\l''}=0$. Thus we have
\begin{eqnarray*}
c_{\theta_0,\l''}&=&(-1)^{l(m_{\theta_0})-l(M_{\theta_0})}\pi_{\theta_0}
\overline{b}_{\theta_0,\l''}+a_{\theta_0,\n_0}(-1)^{l(m_{\n_0})-l(M_{\n_0})}\pi_{\n_0}
\overline{b}_{\n_0,\l''}\\
&&+a_{\theta_0,\r_0}(-1)^{l(m_{\r_0})-l(M_{\r_0})}\pi_{\r_0}
\overline{b}_{\r_0,\l''}\\
&=&\overline{b}_{\theta_0,\l''}-v^{-2}(-\xi)v^{4}+v^{-4}(v^{1}+v^{3})\\
&=&\overline{b}_{\theta_0,\l''}+v^{1}+v^{3}+v^{-1}+v^{-3}.
\end{eqnarray*}
By Lemma 6.3, we can get $b_{\theta_0,\l''}=0$.

(4) When $m\geq8$, we see that $0<\zeta_j<\l''$ if and only if
$j=0,1,2$. When $m=7$, we see that $0<\zeta_j<\l''$ if and only if
$j=0,1$.

One can get $a_{\zeta_2,\l}=0$ for $\l=\n_0, \r_0$, and $\l''$. Thus
we can get $b_{\zeta_2,\l''}=0$ by Lemma 6.3. With the same reason,
we can get $b_{\zeta_1,\l''}=0$, too.

For $\zeta_0$, we have $a_{\zeta_0,\n_0}=v^{-4}$ and
$a_{\zeta_0,\r_0}=a_{\zeta_0,\l''}=0$. Thus we have
\begin{eqnarray*}
c_{\zeta_0,\l''}&=&(-1)^{l(m_{\zeta_0})-l(M_{\zeta_0})}\pi_{\zeta_0}
\overline{b}_{\zeta_0,\l''}+a_{\zeta_0,\n_0}(-1)^{l(m_{\n_0})-l(M_{\n_0})}\pi_{\n_0}
\overline{b}_{\n_0,\l''}\\
&=&(-\xi)\overline{b}_{\zeta_0,\l''}+v^{-4}(-\xi)v^{4}.
\end{eqnarray*}
By Lemma 6.3, we can get $b_{\zeta_0,\l''}=0$.

(5) If $\phi=(\frac{2m}{3}-k)\a+j\b<\l''$ in $\Lambda^+$ with $5\leq
k< \frac{2m}{3}$. Easily, we can have $a_{\phi, \l}=0$ for $\l=\n_0,
\r_0$, and $\l''$. Thus we can get $b_{\phi,\l''}=0$ by Lemma 6.3.

We finish the proof. \hfill$\Box$\

\noindent {\bf Proposition 6.9} Let
$\l''=mx+y=\frac{2m+1}{3}\a+\frac{m+2}{3}\b\in Y$ with $m\geq4$.
Then for $0<\l<\l''$ in $\Lambda^+$, we have
$${b_{\l,\l''}}=\left\{\begin{array}{ll}
{-v^{-2}}, & \ \ \textrm{if} \;  \l=\l''-\a ,\\
{v^{-1}}, & \ \ \textrm{if} \;  \l=\l''-(\a+\b),\\
{v^{-4}}, & \ \ \textrm{if} \; \l=\l''-(2\a+\b),\\
{0} ,   & \ \ \textrm{otherwise} \;  .\end{array} \right.$$\\

\noindent {\bf Proof.} We set
$\r_j=(\frac{2m+1}{3}-1)\a+\frac{m-1}{3}\b+j\b$,
$\n_j=(\frac{2m+1}{3}-2)\a+\frac{m-1}{3}\b+j\b$,
$\theta_j=(\frac{2m+1}{3}-3)\a+\frac{m-4}{3}\b+j\b$,
$\zeta_j=(\frac{2m+1}{3}-4)\a+\frac{m-4}{3}\b+j\b$ in $\Lambda^+$,
with $j\geq0$.

(1) We see that $\r_j<\l''$ if and only if $j=0,1$. Since $\r_1$ is
not in $X$ and $\l''-\r_1=\a$, we have
\begin{eqnarray*}
c_{\r_1,\l''}&=&\overline{b}_{\r_1,\l''}+a_{\r_1,\l''}(-1)^{l(m_{\l''})-l(M_{\l''})}\pi_{\l''}
\overline{b}_{\l'',\l''}\\
&=&\overline{b}_{\r_1,\l''}-v^{-2}.
\end{eqnarray*}
By Lemma 6.3, we can get $b_{\r_1,\l''}=-v^{-2}$.

Since $\r_0\in X$, $\r_1-\r_0=\b$, and $\l''-\r_0=\a+\b$, we get
\begin{eqnarray*}
c_{\r_0,\l''}&=&(-\xi)\overline{b}_{\r_0,\l''}+a_{\r_0,\r_1}
(-1)^{l(m_{\r_1})-l(M_{\r_1})}\pi_{\r_1}
\overline{b}_{\r_1,\l''}\\
&&+a_{\r_0,\l''}(-1)^{l(m_{\l''})-l(M_{\l''})}\pi_{\l''}\\
&=&(-\xi)\overline{b}_{\r_0,\l''}-v^{-2}(-v^{2})-v^{-2}.
\end{eqnarray*}
By Lemma 6.3, we can get $b_{\r_0,\l''}=v^{-1}$.

(2) We see that $\n_j<\l''$ if and only if $j=0,1$. Since
$\n_1\nless \r_0$, $\r_1-\n_1=\a$, and $\l''-\n_1=2\a$, we get
\begin{eqnarray*}
c_{\n_1,\l''}&=&(-1)^{l(m_{\n_1})-l(M_{\n_1})}\pi_{\n_1}
\overline{b}_{\n_1,\l''}+a_{\n_1,\r_1}(-1)^{l(m_{\r_1})-l(M_{\r_1})}\pi_{\r_1}
\overline{b}_{\r_1,\l''}\\
&=&(-1)^{l(m_{\n_1})-l(M_{\n_1})}\pi_{\n_1}
\overline{b}_{\n_1,\l''}-v^{-2}(-v^{2}).
\end{eqnarray*}
By Lemma 6.3, we can get $b_{\n_1,\l''}=0$.

Since $\n_0$ is not in $X$, $\r_0-\n_0=\a$, $\r_1-\n_0=\a+\b$, and
$\l''-\n_0=2\a+\b$, we get
\begin{eqnarray*}
c_{\n_0,\l''}&=&\overline{b}_{\n_0,\l''}+
a_{\n_0,\r_0}(-1)^{l(m_{\r_0})-l(M_{\r_0})}\pi_{\r_0}
\overline{b}_{\r_0,\l''}\\
&&+a_{\n_0,\r_1} (-1)^{l(m_{\r_1})-l(M_{\r_1})}\pi_{\r_1}
\overline{b}_{\r_1,\l''}+a_{\n_0,\l''}(-1)^{l(m_{\l''})-l(M_{\l''})}\pi_{\l''}\\
&=&\overline{b}_{\n_0,\l''}+v^{-2}-v^{-2}+1+v^{-4}.
\end{eqnarray*}
By Lemma 6.3, we can get $b_{\n_0,\l''}=v^{-4}$.

(3) We see that $0<\theta_j<\l''$ if and only if $j=0,1,2$ and
$m\geq5$. Since $\r_1-\theta_2=2\a$ and $\l''-\theta_2=3\a$, we can
get
$c_{\theta_2,\l''}=(-1)^{l(m_{\theta_2})-l(M_{\theta_2})}\pi_{\theta_2}
\overline{b}_{\theta_2,\l''}$. Thus we get $b_{\theta_2,\l''}=0$ by
Lemma 6.3.

Since $\n_0-\theta_1=\a, \r_0-\theta_1=2\a$, $\r_1-\theta_1=2\a+\b$,
and $\l''-\theta_1=3\a+\b$, we have
\begin{eqnarray*}
c_{\theta_1,\l''}&=&(-1)^{l(m_{\theta_1})-l(M_{\theta_1})}\pi_{\theta_1}
\overline{b}_{\theta_1,\l''}+a_{\theta_1,\n_0}(-1)^{l(m_{\n_0})-l(M_{\n_0})}
\pi_{\n_0}\overline{b}_{\n_0,\l''}\\
&&+a_{\theta_1,\r_1} (-1)^{l(m_{\r_1})-l(M_{\r_1})}\pi_{\r_1}
\overline{b}_{\r_1,\l''}\\
&=&(-1)^{l(m_{\theta_1})-l(M_{\theta_1})}\pi_{\theta_1}\overline{b}_{\theta_1,\l''}
-v^{-2}v^{4}+v^{-4}(-v^{2})\\
&=&(-1)^{l(m_{\theta_1})-l(M_{\theta_1})}\pi_{\theta_1}\overline{b}_{\theta_1,\l''}
-v^{2}-v^{-2}.
\end{eqnarray*}
By Lemma 6.3, we can get $b_{\theta_1,\l''}=0$.

Since $\theta_0\in X$, $\n_0-\theta_0=\a+\b, \r_0-\theta_0=2\a+\b$,
$\r_1-\theta_0=2\a+2\b$, and $\l''-\theta_1=3\a+2\b$, we have
\begin{eqnarray*}
c_{\theta_0,\l''}&=&(-\xi)\overline{b}_{\theta_1,\l''}
+a_{\theta_0,\n_0}(-1)^{l(m_{\n_0})-l(M_{\n_0})}
\pi_{\n_0}\overline{b}_{\n_0,\l''}\\
&&+a_{\theta_0,\r_0} (-1)^{l(m_{\r_0})-l(M_{\r_0})}\pi_{\r_0}
\overline{b}_{\r_0,\l''} +a_{\theta_0,\r_1}
(-1)^{l(m_{\r_1})-l(M_{\r_1})}\pi_{\r_1}
\overline{b}_{\r_1,\l''}\\
&=&(-\xi)\overline{b}_{\theta_0,\l''}
-v^{-2}v^{4}+v^{-4}(-\xi)v-v^{-6}(-v^{2})\\
&=&(-\xi)\overline{b}_{\theta_0,\l''} -v^{2}-v^{-2}-v^{-4}+v^{-4}.
\end{eqnarray*}
By Lemma 6.3, we can get $b_{\theta_0,\l''}=0$.

(4) We see that $0<\zeta_j<\l''$ if and only if $j=0,1,2$ and
$m\geq6$. Since $a_{\zeta_2,\l}=0$ for $\l=\r_1$ and $\l''$, we can
get $b_{\zeta_2,\l''}=0$. With the similar reason, we can get
$b_{\zeta_1,\l''}=0$, too.

For $\zeta_0$, we have $\n_0-\zeta_0=2\a+\b$, and $a_{\zeta_0,\l}=0$
for $\l=\r_0,\r_1$ and $\l''$. Thus we get
\begin{eqnarray*}
c_{\zeta_0,\l''}&=&(-1)^{l(m_{\zeta_0})-l(M_{\zeta_0})}\pi_{\zeta_0}
\overline{b}_{\zeta_0,\l''}+a_{\zeta_0,\n_0}(-1)^{l(m_{\n_0})-l(M_{\n_0})}
\pi_{\n_0}\overline{b}_{\n_0,\l''}\\
&=&\overline{b}_{\zeta_0,\l''}+v^{-4}v^{4}.
\end{eqnarray*}
By Lemma 6.3, we can get $b_{\zeta_0,\l''}=0$.

(5) For $\phi=(\frac{2m+1}{3}-k)\a+j\b<\l''$ in $\Lambda^+$, with
$k\geq5$. We have that $a_{\phi,\l}=0$ for $\l=\n_0,\r_0,\r_1$ and
$\l''$. Thus we can get $b_{\phi,\l''}=0$ by Lemma 6.3.

We finish the proof. \hfill$\Box$\

\noindent {\bf Theorem 6.10} For any $u<w$, $u\in c_0$, $w\in c_1$
such that $\mid L(u)\mid=\mid R(u)\mid=1$, we have $\mu(u,w)=0$.

\noindent {\bf Proof.} By Lemma 1.3, we see that $u\leq_L w$ and
$u\leq_R w$ if $\mu(u,w)\neq0$. Then by Lemma 1.4, we get
$L(w)\subseteq L(u)$ and $R(w)\subseteq R(u)$. Since $\mid
L(u)\mid=\mid R(u)\mid=1$, we get $L(w)= L(u)$ and $R(w)= R(u)$.
Without loss of generality, we assume that $L(w)= L(u)=\{r\}$.

Then by Proposition 6.4 (8-9), we see that $\mu(u,w)=\mu(u',w')$,
where $u'=m_{\l}$ and $w'=m_{\l''}$ for some $\l<\l''\in \Lambda^+$
satisfying $\l\in Z$ and $\l''\in X\bigcup Y$. By Propositions
6.6-6.9, we see that
$\mu(u',w')=\mu(m_{\l},m_{\l''})=\textrm{Res}_{v=0}(b_{\lambda,\lambda''})=0$.

Thus the theorem holds.\hfill$\Box$\

\noindent {\bf Theorem 6.11} For any $u<w$, $u\in c_0$, $w\in c_1$
such that $l(w)-l(u)\geq5$, we have $\mu(u,w)=0$.

\noindent {\bf Proof.} By the proof of the above Theorem, we see
that $L(w)\subseteq L(u)$ and $R(w)\subseteq R(u)$ if
$\mu(u,w)\neq0$.

If $\mid L(u)\mid=2$, without loss of generality, we assume that
$L(u)=\{s,t\}$ and $L(w)=\{s\}$. When $L(sw)=\{t\}$, we have
$\mu(u,w)=\mu(ru,rw)$ by using left star operation with respect to
$\{s,r\}$. We see that $rw\in c_1, ru\in c_0$,
$l(rw)-l(ru)=l(w)-l(u)$, and $L(ru)=L(rw)=\{r\}$. When
$L(sw)=\{r\}$, we have $\mu(u,w)=\mu(ru,sw)$ by using left star
operation with respect to $\{s,r\}$. We see that $sw\in c_1, ru\in
c_0$, $l(sw)-l(ru)=l(w)-l(u)-2$, and $L(ru)=L(sw)=\{r\}$. Similarly,
we can discuss the same problem on the right.

If $l(w)-l(u)=5$,  $\mid L(u)\mid=\mid R(u)\mid=2$, and
$\mu(u,w)=\mu(u',w')$  for some $u'\in c_0,w'\in c_1$ satisfying
$\mid L(u')\mid=\mid R(u')\mid=1$ and $l(w')-l(u')=l(w)-l(u)-4=1$.
Then the $u',w'$ must $u'\nless w'$, which implies that
$\mu(u,w)=\mu(u',w')=0$.

For other cases, if $\mid L(u)\mid=2$ or $\mid R(u)\mid=2$, we can
get $\mu(u,w)=\mu(u',w')$  for some $u'\in c_0,w'\in c_1$ satisfying
$u'<w'$, $\mid L(u')\mid=\mid R(u')\mid=1$. By Theorem 6.10, we can
get $\mu(u,w)=\mu(u',w')=0$.

We finish the proof.\hfill$\Box$\

\noindent {\bf Corollary 6.12} The $W$-graph for an affine Weyl
group of type $\widetilde{A}_2$ is locally finite.

\noindent {\bf Proof.} Recall the definition of $W$-graph in [KL1].
Combing Lemma 3.1, Corollary 4.7, Theorems 4.8, 5.12, and 6.11, we
see that, for any given $w\in  W$, there are only finitely many
elements $u\in  W$ such that $\widetilde{\mu}(u,w)\neq0$. Thus we
get the result. \hfill$\Box$\
\providecommand{\bysame}{\leavevmode\hbox
to3em{\hrulefill}\thinspace}
\providecommand{\MR}{\relax\ifhmode\unskip\space\fi MR }
\providecommand{\MRhref}[2]{%
  \href{http://www.ams.org/mathscinet-getitem?mr=#1}{#2}
} \providecommand{\href}[2]{#2}

\end{document}